\title{ The geometrical quantity in damped wave equations on a square}

\author{Pascal H\'ebrard and  Emmanuel Humbert\thanks{ Institut \'Elie
 Cartan, Universit\'e de Nancy 
   1, BP 239  54506 Vandoeuvre-l\`es-Nancy Cedex, FRANCE -
Email: pascal\_hebrard@ds-fr.com and humbert@iecn.u-nancy.fr }}
\date{ November 2003 }
\documentclass[10pt]{article}
\usepackage{amsmath}
\numberwithin{equation}{section}

\textheight 20 cm\usepackage{amssymb,amsmath,epsfig}
\usepackage{graphics}

\newtheorem{theorem}{Theorem}[section]

\newtheorem{proposition}{Proposition}[section]
\newtheorem{definition}{Definition}[section]
\newtheorem{remark}{Remark}[section]

\newtheorem{step}{Step}

\def\L{\Delta}
\def\Om{\Omega}
\def\om{\omega}
\def\l{\lambda}
\def\R{\mathbb{R}}
\def\Z{\mathbb{Z}}
\def\N{\mathbb{N}}
\def\Q{\mathbb{Q}}
\def\p{\partial}

\def\Pr{{\bf Proof~: }}
\def\fpr{\nopagebreak \hfill $\square$}

\begin{document}

\maketitle
\begin{abstract}
The energy in a square membrane $\Om$ subject to constant viscous damping
on a subset $\om \subset \Om$ decays exponentially in time
as soon as $\om$
satisfies a geometrical condition known as the ``Bardos-Lebeau-Rauch'' condition. The rate
$\tau(\om)$ of this decay satisfies $\tau(\om)= 2 \min( -\mu(\om),  g(\om))$
(see Lebeau \cite{leb1}). Here $\mu(\om)$ denotes the spectral abscissa of the
damped wave equation operator and  
$g(\om)$ is a number called 
the geometrical quantity of $\om$ and defined as follows.
A ray in $\Om$ is the trajectory generated by the
free motion of a mass-point in $\Om$ subject to elastic reflections on the
boundary. These reflections obey the law of geometrical optics.
The geometrical quantity $g(\om)$ is then defined as the upper limit (large time
asymptotics) of the average trajectory length.
We give here an algorithm to compute explicitly 
$g(\om)$ when $\om$
is a finite union of squares. 
\end{abstract} 
 
{\bf {\center MSC Numbers: 35L05, 93D15}}
{\bf {\center Key words: Damped wave equation, mathematical billards}}\\

\noindent Let $\Om = [0,1]^2$ be the unit square of $\R^2$ and let $\omega
\subset \Om$ be a
subdomain of $\Om$. We are interested here
in the problem of uniform stabilization of solutions of the following
equation

\begin{equation} \label{eq_pb_lebeau} 
\begin{cases} 
u_{tt}(x,t)-\L u(x,t)+ 2\chi_{\omega} (x)u_{t}(x,t)=0,\; x\in
\stackrel{\circ}\Om, \,t 
> 0,\\ 
u(x,t)=0,\; x\in \p \Om, \,t > 0\\ u(.,0)=u_0 \in H^1_0(\Om),\\ 
u_t(.,0)=u_1 \in L^2(\Om), 
\end{cases} 
\end{equation} 
where $\chi_{\omega}$ is the characteristical function of the subset
$\omega$. This equation has its origin in a physical problem. 
Consider a square membrane $\Om$. We study here the behaviour of a wave in $\Omega$. Let $u(x,t)$
be the vertical position of $x \in \Om$ at time $t>0$. We
assume that we apply on $\omega$ a force proportional to the speed of the
membrane at
$x$. Then, $u$ satisfies equation (\ref{eq_pb_lebeau}). To get 
more information
on this subject, one can refer to \cite{heb}, \cite{rt}, \cite{blr}. With
regard to  the
one-dimensional case, the reader can consult \cite{cz} and
\cite{br}.

\noindent 
We do not deal here with the existence of solutions. We assume that there
exists a solution $u$. Let us define the energy of $u$ at time $t$ by 
$$E(t) = \int_{\Om} {|\nabla u |}^2 + u_t^2 dx.$$
It is well known that, for any $t \geq 0$,

\begin{equation} \label{a_1} 
E(t) \leq C E(0) \hbox{e}^{-2 \tau t },
\end{equation} 
where $C, \tau \geq 0$. As proven by E. Zuazua \cite{zua}, 
this result remains  true  in the semilinear case
when  the  dissipation is in a neighborhood of a subset of the boundary
satisfying the multiplier condition. Let us now define

\begin{definition} 
The exponential rate of decay  $\tau(\omega)$ is defined by 
$$ \tau(\omega) = sup \{\tau \geq 0 \hbox{ s.t. } \exists C> 0 \hbox{ for
  which } (\ref{a_1}) \hbox{ holds } \}.$$
\end{definition}

\noindent Many articles have been devoted to finding bounds for
$\tau(\omega)$. The reader can consult \cite{al}, \cite{leb1} and \cite{sj}.
In the case of a non-constant damping, the reader may see  \cite{cc}.
G. Lebeau proved  in  \cite{leb1} that  
 $$\tau(\om)= 2\min( -\mu(\om), g(\om)).$$
 Here $\mu(\om)$ denotes the spectral abscissa of the
damped wave equation operator and  
$g(\om)$ is a number called the geometrical quantity of $\om$ and defined as follows. 
A ray in $\Om$  is the trajectory generated by the
free motion of a mass-point in $\Om$ subject to elastic reflections in the
boundary. These reflections obey  the law of geometrical optics: the
angle of incidence equals the angle of reflection. If a ray meets a corner
of $\Om$, the reflection will be the limit of the reflection of the rays 
which go
to this corner. It is easy to verify that the ray runs along the same
trajectory before
and after the reflection  but in an opposite direction.

Let $\gamma: \R^+ \rightarrow \Om$ be the parametrization of the ray by
arclength. 
Let 
$$\mathcal{C} = \big\{ \hbox{ rays in } \Om \big\}.$$
Furthermore, let  $\rho_0=(X_0,\alpha)$, where $X_0 \in \Om$ and $\alpha \in [0,2
\pi[$. 
Consider
the ray $\gamma_{\rho_0}$ which starts at $X_0$ in the direction of
the vector  $(\cos(\alpha), \sin(\alpha))$.  
Each ray of $\mathcal{C}$ can be defined in this way. More
precisely, if 
$\Gamma= \Om \times [0,2\pi[ $
then 
$$ \mathcal{C} = \{ \gamma_{\rho_0} \hbox{ s.t. } \rho_0 \in \Gamma \}.$$
In the whole paper  $\gamma \in \mathcal{C}$ will be noted  
$\gamma=[X_0,\alpha]$. 
If $\gamma \in \mathcal{C}$ and if $t > 0$ is a positive real number, 
we write:
$$\gamma_t=\gamma_{/[0,t]}.$$
Let 

$$\mathcal{C}'=\bigcup_{t > 0} \big\{ \gamma_t  \big| \gamma \in
\mathcal{C} \big\} $$
be the set of paths of finite length not necessarily closed.
For $\gamma_t \in \mathcal{C}'$ of length $t>0$, we define:
$$m(\gamma_t)= \frac{1}{t} \int_0^t \mathcal{X}_{\omega}(\gamma_t(s)) ds.$$
The definition of $m$ can be extended  to $\mathcal{C}$. Indeed, if $\gamma
\in \mathcal{C}$, set
$$m(\gamma)=\limsup_{t \to +\infty} m(\gamma_t).$$
It is proven in Section 1.2 that 
\begin{equation} \label{defm}
m(\gamma)=\lim_{t \to +\infty} m(\gamma_t).
\end{equation}
For a ray $\gamma$ belonging to $\mathcal{C}$,
$m(\gamma)$ represents the average time that $\gamma$ spends in $\Om$.

\begin{definition} 
The geometrical quantity $g(\omega)$ is defined by  
$$g(\omega)=  \limsup_{t \to +\infty}
\inf_{\gamma \in \mathcal{C}} m(\gamma_t).$$ 
\end{definition} 

\noindent We are interested here in a precise study of the
geometrical quantity $g$. The first part of
this paper is devoted to studying the rays. We recall some
well known properties of rays. In the second part, another expression
for $g$ is given. Namely, we prove  that 
if  $\om \subset \Om$ whose boundary is a finite union of $C^1$-curves,
then 
$$ g(\omega)=  
\inf_{\gamma \in \mathcal{C}}\limsup_{t \to +\infty}   m(\gamma_t)=
\inf_{\gamma \in \mathcal{C}} m(\gamma).$$

\noindent It is easier to work with this second definition of
$g(\omega)$.
An important application of this theorem is given in the third part: 
we obtain an algorithm which gives an exact computation of $g(\omega)$ when
$\om$ is a finite union of squares. This work has many interests. At first,
Theorem \ref{th_egcrit} gives another definition for $g(\om)$, much more 
easy to manipulate than the original one. Secondly, maximizing the exponential 
rate of decay is interesting from the point of view of physics. 
For these questions,
knowing exactly $g(\om)$ is important. We give some exact computations of
$g(\om)$ with the help of our algorithm at the end of the paper. An
interesting question is then: how can we choose 
the  subset $\omega$ such that 
$g(\omega)$ is maximum? At the moment, this problem is still open. 
The following inequality is always true
 $g(\omega) \leq |\omega|$ (see Section 2) where $|\omega|$ is the area of
 $\omega$. Let us set for $\alpha \in [0,1]$
$$\mathcal{S}(\alpha)= \sup_{|\omega|=\alpha} g(\omega).$$
Some questions then arise naturally

\noindent -  What is the value of $\mathcal{S}(\alpha)$? In particular, is
it true that $\mathcal{S}(\alpha) = \alpha$?

\noindent - Can we find a subset $\omega$ for which $g(\omega)= |\omega|$?

\noindent  - Can we  find an optimal $\om$ (i.e.\  an $\om$ that   maximizes
$g$)  among the domains 
that  satisfy the
multiplier condition?   

\noindent If $|\omega|=0$ or $|\omega|=1$, the answers  are
obvious. However,  if 
$|\omega| \in ]0,1[$, these questions seem to be much more difficult and
are still open.



\section{Basic properties of rays} 
\subsection{Different representations of rays }

In this section, we regard billards in $\Om$ from different point of views
(see for example \cite{tab}). 
Instead of reflecting the trajectory with respect to a side of $\Om$, one
can reflect
$\Om$ with respect to this side. We then obtain a square grid and the initial trajectory
is straightened to a line. This gives  a correspondance between
billard trajectories in $\Om$ and straight lines in the plane equipped with a
square grid.  Two lines in the plane correspond to the same billard
trajectory if they differ by a translation through a vector of the lattice
$2 \Z+2 \Z$. The factor $2$ in $2\Z+2 \Z$ is important. Indeed, two
adjacent squares have an 
opposite orientation in the sense that they are symmetric with respect to
their common side. \\

\begin{figure}[h] 
\centerline{ \psfig{figure=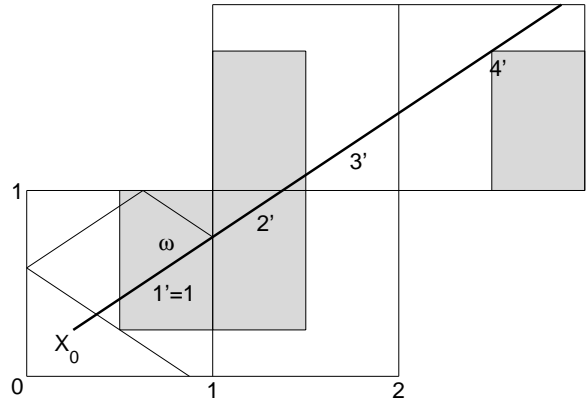,height=6cm}} 
\caption{Point of view of infinite plane}
\label{plane}
\end{figure} 

\noindent Now, consider $T ={[0,2]}^2$ and identify its opposite sides. 
A  trajectory then becomes a geodesic
line in the flat torus.\\

\begin{figure}[h] 
\centerline{ \psfig{figure=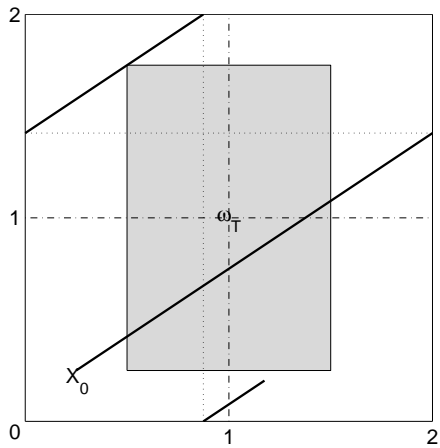,height=6cm}} 
\caption{Point of view of flat torus } 
\label{torus} 
\end{figure}

\noindent Finally, trajectories  in $\Om$ can be seen in three
different ways

\noindent - the definition which consists in reflecting trajectories in
sides of $\Om$

\noindent - the point of view of straight lines in $\R^2$

\noindent - the point of view of geodesic lines in a flat torus $T$.

\noindent Note that  the set $\omega$ must be reflected in the same way
as we did for 
$\Om$. This means that the sets  $\omega$ contained in two adjacent 
squares must be symmetric with respect to their common side (see Figures
\ref{plane} and \ref{torus}).\\

\noindent Let  $t>0$ and $\gamma \in
\mathcal{C}$. The value 
$$\frac{1}{t} \int_0^t \mathcal{X}_{\omega}(\gamma_t(s)) ds$$
does not depend on the point of view we adopt.

\subsection{Open and closed rays} 
 
A ray is said to be closed if it is periodic in the flat torus.  
It is said to be open if it is not closed.
 We have the following characterization. The ray 
 $[X_0=(x_0,y_0),\alpha]$ 
is open if and only if  $\cos \alpha$ and  $\sin \alpha$ are independent
over $\Z$, i.e. if and only if $\tan \alpha \in \R-\Q$. 
Indeed, for example using the point of view  of a flat torus, the trajectory is
periodic if and only if there exist  $0 \leq t_1 <t_2$ and  $n,m \in \N$  
such that
$$ \left\{ \begin{array}{ll} x_0 + t_1 \cos \alpha = 
x_0 + t_2 \cos \alpha +2n\\ y_0 + t_1 \sin \alpha = y_0 + t_2 \sin 
\alpha + 2m 
\end{array} \right. .$$ 
Hence $2m(t_1-t_2) \cos \alpha - 
2n(t_1-t_2) \sin \alpha = 0$.
\noindent Let us now consider open rays. They have the following properties
(see for example \cite{cfm} p. 172).
\begin{proposition}\label{prop1} 
Let  $\gamma = [X_0, \alpha]$ be an open ray. Then,
$\{\gamma_t, \, t\geq 0\}$ is dense in the torus  $T$ and then
in $\Om$. In addition, if $\om$ is  quarrable (i.e. $\chi_\om$
Riemann-integrable on $\Om$) 
 then  $m(\gamma_t) 
\xrightarrow[t \to + \infty]{}|\om|$. More precisely, let $\alpha$ 
be such that $\tan \alpha \in \R - \Q$ then
$$\forall \epsilon > 0, \, \exists t_0\, |\, t\geq t_0,
\forall X_0 \in \Om \; : \; \text{if} \; \gamma_{\rho_0} = 
[X_0,\alpha] \; \text{then} \; \left|m(\gamma_t)-|\om|\right| < \epsilon.$$ 
In other words,  $t_0$ is independent of  $X_0$. 
\end{proposition} 
Here, $|\om|$ stands for the area of $\om$. 
Consider now closed rays. They satisfy the following properties:

\begin{proposition} \label{prop} 
Let $\gamma=[X_0, \alpha]$ be a closed ray.
Then, there exist  $p,q \in \Z$ relatively prime
such that: 
$$ \cos \alpha = \frac{p}{\sqrt{p^2+q^2}} \quad \text{and} 
\quad \sin \alpha = \frac{q}{\sqrt{p^2+q^2}}.$$ 
The period $L$ of the trajectory is $L=2\sqrt{p^2+q^2}$. We have
$$ m(\gamma_t) \xrightarrow[t \to + \infty]{} \frac{1}{L}\int_0^L 
\chi_\om(\gamma(s))ds. $$ 
\end{proposition} 

\begin{remark} \label{rem1}
The assertion (\ref{defm}) follows immediately from Propositions \ref{prop1}
and \ref{prop}. Moreover, let $\gamma \in \mathcal{C}$. If $\gamma$ is
open, we have $m(\gamma)=|\om|$ and if $\gamma$ is periodic of period $L$,
we have $m(\gamma)=  \frac{1}{L}\int_0^L 
\chi_{\om}(\gamma(s))ds $.
\end{remark}

\noindent As one can check, an immediate consequence of Proposition 
\ref{prop} is the following result 
\begin{proposition}\label{th_majgeom} 
Let $\om \subset \Om$ be  Riemann-integrable, then 
$$g(\om) \leq |\om|.$$ 
\end{proposition}

\noindent Closed rays can easily be described. Let us adopt 
the point of view of a flat torus. Let us  consider a ray  $\gamma=[X_0,\alpha]$  with $X_0 \in T$ and $\cos 
\alpha = p/ \sqrt{p^2+q^2},\, \sin \alpha = q/ \sqrt{p^2+q^2}$, 
$p$ and $q$ relatively prime. We
assume that  $p \geq 0$ and $q \geq 0$  (i.e 
$\alpha \in [0,\pi/2]$). From remark
\ref{rem_points_dep} below,  the study can be restricted to such rays.
\begin{theorem}\label{th_rayfermedansT} 
The ray $\gamma$ is an union of $p+q$ parallel segments directed by 
 $\vec{u} = \cos \alpha \, 
\vec{i} + \sin \alpha \, \vec {j}$. Among these segments, $p$ of them start
at    
$A_k$, $1 \leq k \leq p$, and the $q$  others  start at  
$B_k$, $1 \leq k \leq q$ with $$ A_k \, \left(0, 
\big(y_0+q/p(2k - x_0)- 2E(\frac{1}{2}(y_0+q/p(2k - x_0))\big)\right)$$ 
and 
$$B_k \, 
 \left( 
\big(x_0+p/q(2k - y_0)- 2E(\frac{1}{2}(x_0+p/q(2k - y_0))\big) ,0\right)).$$ 
In addition, the distance between two neighbouring segments is constant and
equal to $\delta = 4/L$, where 
$L = 2\sqrt{p^2+q^2}$ is the period. 
\end{theorem} 
 
\begin{figure}[h] 
\centerline{ \psfig{figure=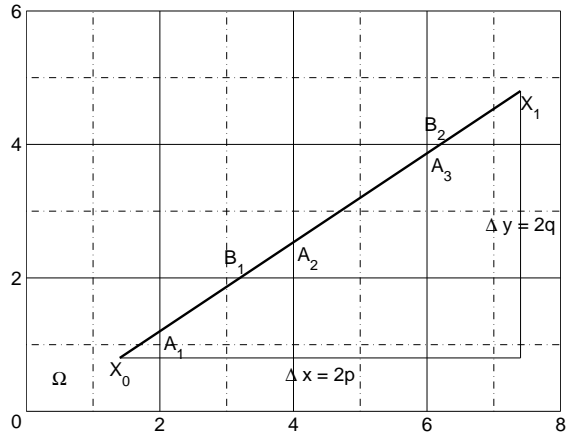,height=6cm}} 
\caption{Closed ray for $p=3$ and $q=2$ in the infinite plane} 
\label{fig_ray_ferme_plan} 
\end{figure} 
\begin{figure}[h] 
\centerline{ \psfig{figure=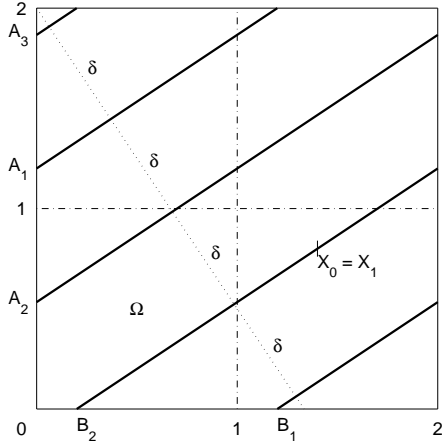,height=6cm}} 
\caption{Closed ray $p=3$ and $q=2$ in the flat torus} 
\label{fig_ray_ferme_tore} 
\end{figure} 
 
\noindent \Pr In the  infinite plane $\R^2$, 
the ray is defined by $x(t) = x_0 +2t p/L$ 
and $y(t) = y_0 + 2t q/L$, where $L = 2\sqrt{p^2+q^2}$ is the period. 
Let us restrict the study to the segment $[X_0 X_1]$ of length $L$. In $]X_0
X_1[$, the ray intersects $p$ 
vertical straight lines defined by the  equations  
$x=2k,\, k = 1, \cdots,p$ and  $q$ horizontal straight lines defined by the
equation 
$y=2k,\, k= 1, \cdots ,q$. This gives $p+q$ segments starting at points with coordinates 
$$ \big(0, 
y_0+q/p(2k - x_0)\big)$$ and 
$$\big( 
x_0+p/q(2k - y_0),0\big).$$
Coming back to the point of view of the torus, this gives  $p+q$ 
segments starting at $A_k$ and 
$B_k$.\\ 
\noindent Since $p$ and  $q$ are relatively prime, it is well known that 
  $\{k q \, [p],\, 
1 \leq k \leq p\} = \{0,1,\dots,p-1\}$. Hence, if $A_{k_1}$ and $A_{k_2}$
are two   
neighbouring points, $A_{k_1}A_{k_2} = 2/p$. As a consequence, the distance
between the straight lines $(A_{k_1},\vec{u})$ and 
$(A_{k_2},\vec{u})$ is $$\delta = \left| 
\overrightarrow{A_{k_1}A_{k_2}}.\vec{u}^\bot \right| = \frac{2}{p} 
\, \vec{j}.\frac{-q \vec{i} + 
p\vec{j}}{\sqrt{p^2+q^2}}=\frac{4}{L}. $$ 
In the same way, if $B_{k_1}$ 
and $B_{k_2}$ are two neighbouring points $B_{k_1}B_{k_2} = 2/q$, the
distance between $(B_{k_1},\vec{u})$ and 
$(B_{k_2},\vec{u})$ is $$\delta = \left| 
\overrightarrow{B_{k_1}B_{k_2}}.\vec{u}^\bot \right| = \frac{2}{q} 
\, \vec{i}.\frac{q \vec{i} - 
p\vec{j}}{\sqrt{p^2+q^2}}=\frac{4}{L}.$$ 
 
\noindent Let $k_1$, $k_2$ be  such that 
$$dist (A_{k_1},0)= \min_k dist(A_k,0)\hbox{ and } dist (B_{k_2},0)= \min_k
dist(B_k,0).$$
It remains to prove that the distance between the two segments starting at
$A_{k_1}$ and $B_{k_2}$ is $\delta = 4/L$. Since the distance between two
points  $A_k$ is $2/p$, there exists $n \in \N$ such that the 
ordinate $y^*$ of $A_{k_1}$ satisfies 
\begin{equation*} \begin{split} 
2n &\leq y_0+q/p(2k_1 - x_0)- 2E(\frac{1}{2}(y_0+q/p(2k_1 - x_0)) < 2n + 2/p\\ 
   &\; \Longrightarrow \; y^*=  y_0+q/p(2k_1 - x_0)- 2E(\frac{1}{2}(y_0+q/p(2k_1 - x_0))  -2n <2/p. 
\end{split} \end{equation*} 
In the same way, there exists $m \in \N$ such that the abscissa $x^*$ 
of $B_{k_2}$ satisfies 
\begin{equation*} \begin{split} 
2m &\leq x_0+p/q(2k_2 - y_0)- 2E(\frac{1}{2}(x_0+p/q(2k_2 - y_0))  < 2m + 2/q \\ 
   &\; \Longrightarrow \; x^*= x_0+p/q(2k_2 - y_0)- 2E(\frac{1}{2}(x_0+p/q(2k_2 - y_0))  - 2m <2/q. 
\end{split} \end{equation*} 
The distance $\delta$ between the two segments verifies:
$$\delta = \left| \overrightarrow{A_{k_1}B_{k_2}}.\vec{u}^\bot 
\right| = \frac{|q x^* + p y^*|}{\sqrt{p^2+q^2}}.$$ 
An easy computation shows that $|q x^* + p y^*|$ is an even integer
number. 
Consequently, $\delta$ can be written as $4k/L$, where $k$ is an integer. 
Since $0\leq x^* < 2/q$ and  $0\leq y^* < 2/p$, we get $\delta < 8/L$, and  
hence, $\delta = 4/L$. \fpr 
 
\begin{remark} \label{rem12}
If a quarrable domain  $\om$ is such that $g(\om) = |\om|$ then 
for all $\alpha \in [0,\pi/2]$ and for almost all $X_0 \in T$, 
$m(\gamma)= |\om|$, where $\gamma 
= [X_0, \alpha]$. 
\end{remark} 
Indeed, consider an open ray $\gamma$. Then 
$m(\gamma)=  |\om|$. 
Let $X_0 \in \Om$. Assume that a ray $\gamma$ starting at $X_0$ 
is closed of length $L$. 
Let $M$ and $N$ be two points of $\gamma$ such that  $\overrightarrow{MN}$ is
orthogonal to  $\vec{u}$ and $MN = 
\delta$. Let also $P(s) = s M +(1-s)N$ and $\gamma^s=[P(s),\alpha]$, 
then $$ g(\om) \leq \inf_{s \in [0,1]} m(\gamma_L^s) \leq 
\int_0^1 m(\gamma_L^s) ds \leq |\om|.$$ 
If $g(\om)=|\om|$, then $ m(\gamma_L^s) = |\om|$ almost everywhere. 
If $\partial \om$  is  a finite
union of curves of class 
$C^1$, then the number of discontinuity points of $ m(\gamma_L^s)$ is
finite. It is null if $\partial \om$ does not possess any segments. 
\noindent The remark above immediately follows.

\section{An equivalent definition for $g$} 
\subsection{ The geometrical quantity $g'$}
 We define, for any set $\om$:

\begin{definition} The geometrical quantity $g'(\om)$ is defined by 
$$g'(\om)= \inf_{\gamma \in \mathcal{C}} m(\gamma). $$
\end{definition} 

\noindent 
As easily seen, $g'$ is easier to study than $g$. As shown in
Section 3, $g'(\om)$ can be computed explicitly for a certain class of domains
$\om$. Together with Theorem \ref{th_egcrit}, this gives an algorithm for
computing $g(\om)$.

\begin{theorem} \label{th_egcrit} 
Let $\om$ be a closed set whose boundary is a finite union of
curves of class $C^1$. Then 
$$ g(\om) = g'(\om).$$
\end{theorem}

\noindent As a first remark, the same argument as in the proof of
Proposition
\ref{th_majgeom} shows that 

\begin{eqnarray} \label{b1}
0 \leq  g'(\om) \leq 
|\omega|.
\end{eqnarray}

\noindent The proof of Theorem \ref{th_egcrit} is given in the appendix. It is
easy to see that 
$g(\om) \leq g'(\om)$. Inequality $g'(\om) \leq g(\om)$ is much more
difficult to obtain. To prove Theorem \ref{th_egcrit}, we consider a
sequence of rays 
$\gamma^n=[x_n, \theta_n]$  and a sequence of real numbers  $t_n \to
+\infty$ for which 
$m(\gamma^n_{t_n}) \to g(\om)$. After choosing  a subsequence, there exists $(x,
\theta) \in \Omega \times [0,2 \pi[$ such that 
$\lim_n x_n =x$ and $\lim_n \theta_n = \theta$. We will now show that 
$\lim_n m(\gamma^n_{t_n}) - m([x,\theta]_{t_n})=0$. The conclusion
then follows. The difficulty in this proof is that the
function $(x,\theta)  \to m([x,\theta]_t)$ is not continuous.  
A direct consequence of Theorem
\ref{th_egcrit} is that it suffices to consider closed rays in the
explicit computation of $g(\om)$. Namely,
let  $\mathcal{C}_c$ be the set of closed rays, then we have: 
\begin{equation} 
 g(\om) = \inf_{\gamma \in \mathcal{C}_c} m(\gamma).
\end{equation}


\section{Explicit computation of  $g$} 
We prove in this section that for a particular class of domains
$\om$, 
 the properties of the geometrical quantity  $g$ we obtain above
allow  to compute explicitly $g(\om)$. Let $N \in \N^*$, and let $\om
\subset \Om$ be a finite union of squares  
$(C_{i,j})_{1 \leq i,j \leq N}$, with $$C_{i,j} = 
\left[ \frac{i-1}{N},\frac{i}{N} \right] \times \left[ 
\frac{j-1}{N},\frac{j}{N} \right].$$ 
Obviously, $\om$ is Riemann-measurable, closed and its boundary is a finite
union of $C^1$ curves.
 
\subsection{Influence of $\alpha$} 
The first result we obtain is the following 
\begin{theorem} \label{th_influencealpha} 
Let $\gamma = [X_0,\alpha]$ be a closed ray  such that 
$\tan \alpha = q/p$, with $p$ and $q$ relatively prime and $p + q 
>2N$, then 
$$\left| \, m(\gamma)- 
|\om| \, \right| \leq \frac{N^2}{p q} \min \{ |\om|, 1-|\om| \} 
\leq \frac{N^2}{2p q}.$$ 
\end{theorem} 
\noindent \Pr We adopt the point of view of flat torus $T=[0, 2]^2$. 
At first, let 
$K$ be one of the $C_{i,j}$ i.e. 
 $$K = 
C_{i_0,j_0} = \left[ \frac{i_0-1}{N}, \frac{i_0}{N} \right] \times 
\left[ \frac{j_0-1}{N}, \frac{j_0}{N} \right] = ABCD.$$
The period of $\gamma$ is   $L=2\sqrt{p^2+q^2}$. Let $\delta = 
4/L$. Then  $\cos \alpha = p\delta/2$ and  $\sin \alpha = 
q\delta/2$. By Theorem
\ref{th_rayfermedansT}, the ray $\gamma$ is constituted by 
$p+q$ segments directed by  $\vec{u}$. The distance between two
neighbouring segments is $\delta$.\\ 

\begin{figure}[h] 
\centerline{\psfig{figure=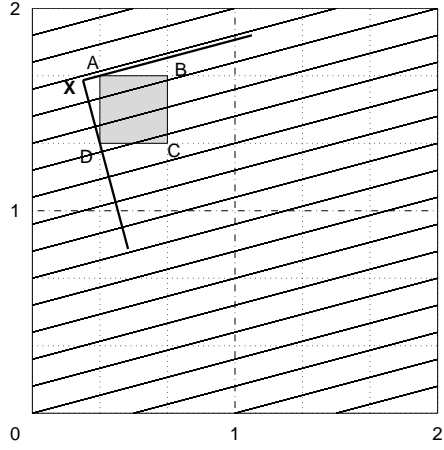,height=6cm}} 
\caption{Ray $\gamma = [0, \alpha]$ with  $\tan \alpha = 
4/15$} \label{fig_repere_oblique} 
\end{figure}

\noindent Let us define  the orthonormal frame 
$\mathcal{R}_u = (X,\overrightarrow{I_u},\overrightarrow{J_u})$, 
where  $X$ is such that $\overrightarrow{XD}.\vec{u} = 0$,
$\overrightarrow{XD}.\overrightarrow{XA} = 0$, 
$\overrightarrow{I_u} = 
\overrightarrow{XD}/||\overrightarrow{XD}||$, and  
$\overrightarrow{J_u} = \overrightarrow{XA}/||\overrightarrow{XA}||
= \vec{u}$ (see Figure \ref{fig_repere_oblique}).\\ 
 
\begin{figure}[h] 
\centerline{\psfig{figure=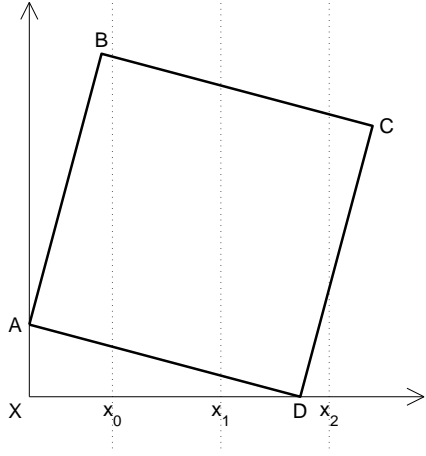,height=6cm}} 
\caption{Vertical segments in the square $ABCD$} 
\label{fig_carre_oblique} 
\end{figure}

\noindent In this frame, the coordinates of $A, B, C$ and $D$ satisfy $$ A\,(0, \frac{\sin \alpha}{N}); \quad B\, 
(\frac{\sin \alpha}{N}, \frac{\sin \alpha + \cos \alpha}{N}); 
\quad C\, (\frac{\sin \alpha + \cos \alpha}{N}, \frac{\cos 
\alpha}{N}); \quad D\,( \frac{\cos \alpha}{N}, 0).$$ 
The ray  
$\gamma$ is constituted of $p+q$ vertical segments. $R$ of them  
meet $K$ (see  Figure 
\ref{fig_carre_oblique}). We note  $x_0, x_1, \dots x_{R-1}$ their  
abscissa $$\left\{ 
\begin{array}{ll} x_0 = r \in [0,\delta] \\ x_i = x_0 + i \delta, 
\quad 0 \leq i \leq R-1 
\end{array} \right. .$$ 
Since $\frac{\sin \alpha + \cos \alpha}{N} =  \frac{p + 
q}{2N}\delta > \delta$, there exists at least one segment in  
$K$ and  $R \geq 1$.\\ 
 
\noindent Let  $$ 
\begin{array}{lllll} f_1(x) & = \cot \alpha \; x + \dfrac{\sin 
\alpha}{N} =  \cot \alpha \, (x - \dfrac{\sin \alpha}{N}) + 
\dfrac{\sin \alpha + \cos \alpha}{N} \; \text{ line } (AB)\\ 
f_2(x) & = -\tan \alpha \, (x - \dfrac{\sin \alpha}{N}) + 
\dfrac{\sin \alpha + \cos \alpha}{N}\\ & =  -\tan \alpha \, (x - 
\dfrac{\sin \alpha + \cos \alpha}{N}) + \dfrac{\cos \alpha}{N} \; 
\text{ line } (BC)\\ g_1(x) & = -\tan \alpha \; x + \dfrac{\sin 
\alpha}{N} =  -\tan \alpha \, (x - \dfrac{\cos \alpha}{N}) \; 
\text{ line } (AD)\\ g_2(x) & = \cot \alpha \, (x - \dfrac{\cos 
\alpha}{N}) =  \cot \alpha \, (x - \dfrac{\sin \alpha + \cos 
\alpha}{N}) + \dfrac{\cos \alpha}{N} \; \text{ line } (CD) . 
\end{array}$$ and $$f(x) = \left\{ 
\begin{array}{llll} \dfrac{\sin \alpha}{N} & \text{if } x \leq 0 
\\ f_1(x) & \text{if } 0 \leq x \leq \frac{\sin \alpha}{N}\\ 
f_2(x) & \text{if } \frac{\sin \alpha}{N} \leq x \leq \frac{\sin 
\alpha + \cos \alpha}{N}\\ \dfrac{\cos \alpha}{N} & \text{if } x 
\geq \frac{\sin \alpha + \cos \alpha}{N} \end{array} \right. 
\text{and } g(x) = \left\{ 
\begin{array}{llll} \dfrac{\cos \alpha}{N} & \text{if } x \leq 0 
\\ g_1(x) & \text{if } 0 \leq x \leq \frac{\cos \alpha}{N}\\ 
g_2(x) & \text{if } \frac{\cos \alpha}{N} \leq x \leq \frac{\sin 
\alpha + \cos \alpha}{N}\\ \dfrac{\cos \alpha}{N} & \text{if } x 
\geq \frac{\sin \alpha + \cos \alpha}{N} \end{array} \right. .$$ 
Then $$ \int_\R f - g = \int_0^{\frac{\sin \alpha + \cos 
\alpha}{N}} (f(x) - g(x))dx = |K| = \frac{1}{N^2} $$ and $$ 
\frac{1}{L} \int_0^L \chi_{K}(\gamma(s))ds = \frac{\delta}{4} 
\sum_{i=0}^{R-1} f(x_i) - g(x_i).$$ Let  $k$ be the number of vertical
segments 
which meet $[AB]$, i.e. the $k$ segments such that
$$ r + (k-1)\delta \leq \frac{\sin \alpha}{N} \leq r + k 
\delta.$$ If $1 \leq k \leq R-1$, if at least one segment meets $[AB]$ and
if at
least one segment meets $[BC]$ then, as one can see 
\begin{align} 
\nonumber 
\frac{1}{L} \sum_{i=0}^{R-1} f(x_i) &= \frac{\delta}{4} \sum_{i=0}^{k-1} 
f_1(x_i) + \frac{\delta}{4} \sum_{i=k}^{R-1} f_2(x_i)\\ \nonumber &= \frac{1}{4} \sum_{i=0}^{k-1} 
\int_{x_i-\delta/2}^{x_i+\delta/2} f_1(x)dx + \frac{1}{4} \sum_{i=k}^{R-1} 
\int_{x_i-\delta/2}^{x_i+\delta/2} f_2(x)dx\\ &= 
\frac{1}{4} \int_{r-\delta/2}^{r+(k-1/2)\delta} f_1(x)dx + 
\frac{1}{4} \int_{r+(k-1/2)\delta}^{r+(R-1/2)\delta} f_2(x)dx . 
\label{eq_sommef} 
\end{align} 
If  $k=0$ or $k=R$, 
equation (\ref{eq_sommef}) remains true. Hence, in all cases 
\begin{equation*} 
\begin{split} 
\frac{1}{L} \sum_{i=0}^{R-1} f(x_i) &= \frac{1}{4} \int_0^{(p+q)\delta/2N} 
f(x)dx - \frac{1}{4} \int_0^{r-\delta/2} f_1(x)dx\\ &\qquad + 
\frac{1}{4} \int_{q\delta/2N}^{r+(k-1/2)\delta} (f_1(x) - f_2(x))dx + 
\frac{1}{4} \int_{(p+q)\delta/2N}^{r+(R-1/2)\delta} f_2(x)dx. 
\end{split} 
\end{equation*} 
In the same way, let $l$ be the number of segments which meet $[AD]$ $$r + (l-1)\delta \leq \frac{\cos 
\alpha}{N} \leq r + l \delta .$$ 
We obtain that, in all cases
\begin{equation*} 
\begin{split} 
\frac{1}{L} \sum_{i=0}^{R-1} g(x_i) &= 
\frac{1}{4} \int_{r-\delta/2}^{r+(l-1/2)\delta} g_1(x)dx + 
\frac{1}{4} \int_{r+(l-1/2)\delta}^{r+(R-1/2)\delta} g_2(x)dx\\ &= 
\frac{1}{4} \int_0^{(p+q)\delta/2N} g(x)dx - \frac{1}{4} \int_0^{r-\delta/2} g_1(x)dx\\ 
&\qquad + \frac{1}{4} \int_{p\delta/2N}^{r+(l-1/2)\delta} (g_1(x) - g_2(x))dx + 
\frac{1}{4} \int_{(p+q)\delta/2N}^{r+(R-1/2)\delta} g_2(x)dx . 
\end{split} 
\end{equation*} 
Finally 
\begin{equation*} 
\begin{split} 
\frac{1}{L} \int_0^L \chi_K(\gamma(s))ds -\frac{1}{4N^2} &= -\frac{1}{4} \int_0^{r-\delta/2} f_1 -g_1 + 
\frac{1}{4} \int_{q\delta/2N}^{r+(k-1/2)\delta} f_1 - f_2 \\ &\qquad+ 
\frac{1}{4} \int_{p\delta/2N}^{r+(l-1/2)\delta} g_1 - g_2 + 
\frac{1}{4} \int_{(p+q)\delta/2N}^{r+(R-1/2)\delta} f_2 - g_2\\ &= \frac{\tan 
\alpha + \cot \alpha}{8}\left[ -(r-\delta/2)^2 + (r+(k-1/2)\delta 
- q\delta/2N)^2 \right.\\ &\qquad \left.+ (r+(l-1/2)\delta - 
p\delta/2N)^2 - (r+(R-1/2)\delta - (p+q)\delta/2N)^2 \right] . 
\end{split} 
\end{equation*} 
However, in this last equality, each  of the four square terms is less than
$\delta/2$  
and $\tan \alpha + \cot \alpha = \dfrac{1}{\sin \alpha \, \cos 
\alpha} = \dfrac{4}{pq \delta^2} = \dfrac{L^2}{4pq}$ hence 
\begin{equation}\label{eq_encadcarre} 
\left| \frac{1}{L} \int_0^L \chi_K(\gamma(s))ds - \frac{1}{4N^2}\right| \leq \frac{1}{4pq} . 
\end{equation} 
 
\noindent We now consider:
 $$ \om_T = 
\bigcup_{(i,j) \in E} C_{i,j},$$ where $E$ is a subset of  
$\{1,2,\dots,2N\}^2$ and  $\text{card}(E) = N^2|\om_T| = 4 N^2 
|\om|$~; then $$m(\gamma) = \frac{1}{L} \int_0^L 
\chi_{\om_T}(\gamma(s))ds = \sum_{(i,j) \in E} \frac{1}{L} 
\int_0^L \chi_{C_{i,j}}(\gamma(s))ds.$$ Hence $$m(\gamma) - |\om| 
= \sum_{(i,j) \in E} \left( \frac{1}{L} \int_0^L 
\chi_{C_{i,j}}(\gamma(s))ds - \frac{1}{4N^2} \right).$$ This leads to~$$ 
\left| \, m(\gamma) - |\om| \, \right| \leq 
\frac{\text{card}(E)}{4pq} = \frac{N^2}{pq}|\om|.$$ If the area of  
$\om$ is greater than 1/2, let  $\om'$ be  the adherence of 
$\Om-\om$ which is also  a finite union of squares 
$C_{i,j}$, Hence$$ \left|\frac{1}{L} \int_0^L 
\chi_{\om'}(\gamma(s))ds - (1-|\om|) \, \right| \leq 
\frac{N^2}{pq}(1-|\om|)$$ or$$\frac{1}{L} \int_0^L 
\chi_{\om'}(\gamma(s))ds = 1 - \frac{1}{L} \int_0^L 
\chi_{\om}(\gamma(s))ds.$$ Finally, we get  $\left| \, m(\gamma) - 
|\om| \, \right| \leq \dfrac{N^2}{pq}(1-|\om|)$. \fpr

\subsection{Representation of  $\om$} 
\begin{figure}[h] 
\centerline{\psfig{figure=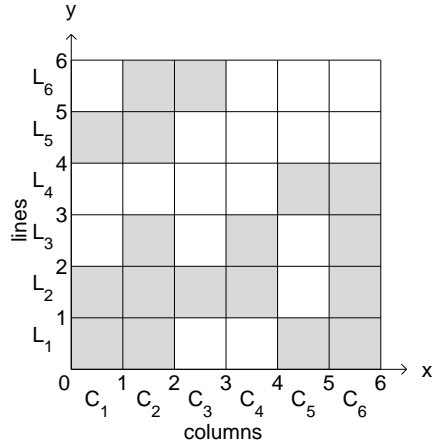,height=6cm}} 
 \caption{A domain $\om$ for  $N = 6$} \label{fig_codage} 
\end{figure}

For convenience, consider the
frame 
 $(O,\vec{I},\vec{J})$ with $\vec{I} = \vec{i}/(2N)$ 
and $\vec{J} = \vec{j}/(2N)$. As a consequence, $\Om$ is now the square 
 $[0,N]^2$, and $T$ is the torus $[0,2N]^2$. The length of closed rays
 defined by  $(p,q)$ ($p,q$ relatively prime) 
 such that  $\cos \alpha = 
p/\sqrt{p^2+q^2}$ and $\sin \alpha = q/\sqrt{p^2+q^2}$ must be modified. 
The period is now
$$L = 2N \sqrt{p^2+q^2}$$ and $C_{i,j}$ is such that: 
$$ C_{i,j} = [i-1,i] \times 
[j-1,j].$$  Let us define the matrix $M$  $N \times N$ such that  
$M_{i,j}=1$ if $ C_{j,i} \subset \om$, and $M_{i,j}=0$ if $ C_{j,i}
\not\subset \om$. With 
this notation, the matrix $M$ is a representation of $\om$. As an
example, if  $\om$ is as in Figure \ref{fig_codage} then 
$$ M = \left( \begin{array}{cccccc} 
 1 & 1 & 0 & 0 & 1 & 1 \\ 
 1 & 1 & 1 & 1 & 0 & 1 \\ 
 0 & 1 & 0 & 1 & 0 & 1 \\ 
 0 & 0 & 0 & 0 & 1 & 1 \\ 
 1 & 1 & 0 & 0 & 0 & 0 \\ 
 0 & 1 & 1 & 0 & 0 & 0 \\ 
\end{array} \right). $$ 
 
\noindent Let us also define the $2N \times 2N$ matrix $A$ such that $A_{i,j}=1$ if
$ C_{j,i} \subset \om_T$, $A_{i,j}= 0$ if $ C_{j,i} \not\subset \om_T$, i.e. 
\begin{equation} \label{eq_defA} 
M_{i,j} = A_{i,j} = A_{2N+1-i,j} = A_{i,2N+1-j} = A_{2N+1-i,2N+1-j}. 
\end{equation}

\subsection{Influent rays} 
Let $\gamma$ be a horizontal ray starting at $X_0=(x_0,y_0)$. Then, 
$$m(\gamma)= \frac{1}{N} \sum_{j=1}^{N} M_{i,j}.$$
Assume that  $y_0=k$ is an integer. Since each $C_{i,j}$ is closed, $$m(\gamma)
= \frac{1}{N} \sum_{j=1}^{N} (M_{k-1,j} + 
M_{k,j}).$$ Obviously, this ray does not realize the infimum in the
definition of $g$. In the same way, for a vertical ray, 
only the numbers $$\frac{1}{N} \sum_{i=1}^{N} M_{i,j}, \quad 1 
\leq i \leq N.$$ must be considered. Now, let us deal with  oblique rays
and 
fix $\alpha$.
Let us denote  $\tan 
\alpha$ by $q/p$, where $p$ and $q$ are relatively prime integers. 
\begin{remark}\label{rem_points_dep} 
One has to consider only the rays defined by  $\gamma =  
[X_0,\alpha]$ with $X_0 \in T$ and  $\alpha \in [0,\pi/2]$, i.e. $p 
\geq 0$ and $q \geq 0$. 
\end{remark} 
 
\begin{figure}[h] 
\centerline{\psfig{figure=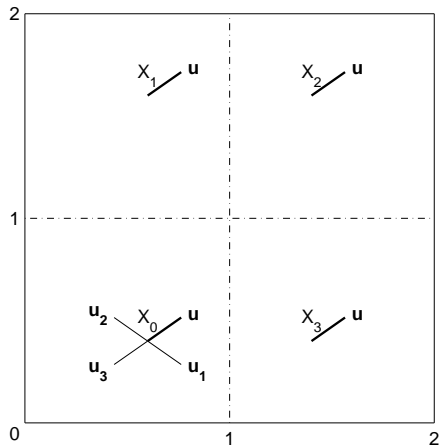,height=6cm}} 
 \caption{Starting points with $\alpha \in [0,\pi/2]$} 
 \label{fig_points_dep} 
\end{figure} 
 
\noindent Indeed, keeping the notations of Figure \ref{fig_points_dep}, the
ray 
 $\gamma^1=[X_0,\overrightarrow{u_1}]$ is 
equivalent to  $\tilde{\gamma}^1=[X_1,\vec{u}]$. In the same way,  
the rays $\gamma^2=[X_0,\overrightarrow{u_2}] $ and
$\gamma_3=[X_0,\overrightarrow{u_3}] $   
are equivalent to
 $\tilde{\gamma}^2=[X_2,\vec{u}] $ and $\tilde{\gamma}^3=[X_3,\vec{u}] $.\\ 
 
\noindent Let us now study the 
influence of the starting point  of
 rays. We have the following result 
\begin{proposition} \label{th_points_nonentiers} 
Among oblique closed rays of angle $\alpha$, the ray which spends the least
time in $\om$ in average is a ray which meets a point with integer
coordinates .
\end{proposition} 
\noindent \Pr Let $B=(x_0,y_0)$ be a point of $T$ such that the ray 
$\gamma=[B,\alpha]$ does not meet  a
point with integer coordinates. 
Let us adopt the point of view of infinite plane and 
let  $B_0, B_1, B_2, \dots$ be the  intersection points of $\gamma$ with the lattice 
$\Z+\Z$.   A direct consequence of
Theorem  \ref{th_rayfermedansT} is that $B_0B_{2N(p+q)} = L 
=2N\sqrt{p^ 2+q^2}$. $\vec{u} = (p \vec{i} + q 
\vec{j})/\sqrt{p^2+q^2}$ is the direction of the ray. Let $\vec{n} 
= (-q \vec{i} + p \vec{j})/\sqrt{p^2+q^2}$ be such that 
$(\vec{u},\vec{n})$ is a direct orthonormal frame. For a point $P=(k,l)$
with integer coordinates, let us  compute the algebraic distance to the ray
$\gamma$  
$$d(P,\gamma) = \overrightarrow{BP}.\vec{n} = \frac{pl-qk 
- (py_0-qx_0)} {\sqrt{p^2+q^2}}.$$ There exists a point $A$  with integer
coordinates and which verifies:   
$$d(A,\gamma)= 
\min \{ d(P,\gamma) | P \text{ has integer coordinates and } d(P,\gamma)>0
\}. $$
 
\begin{figure}[h] 
\centerline{\psfig{figure=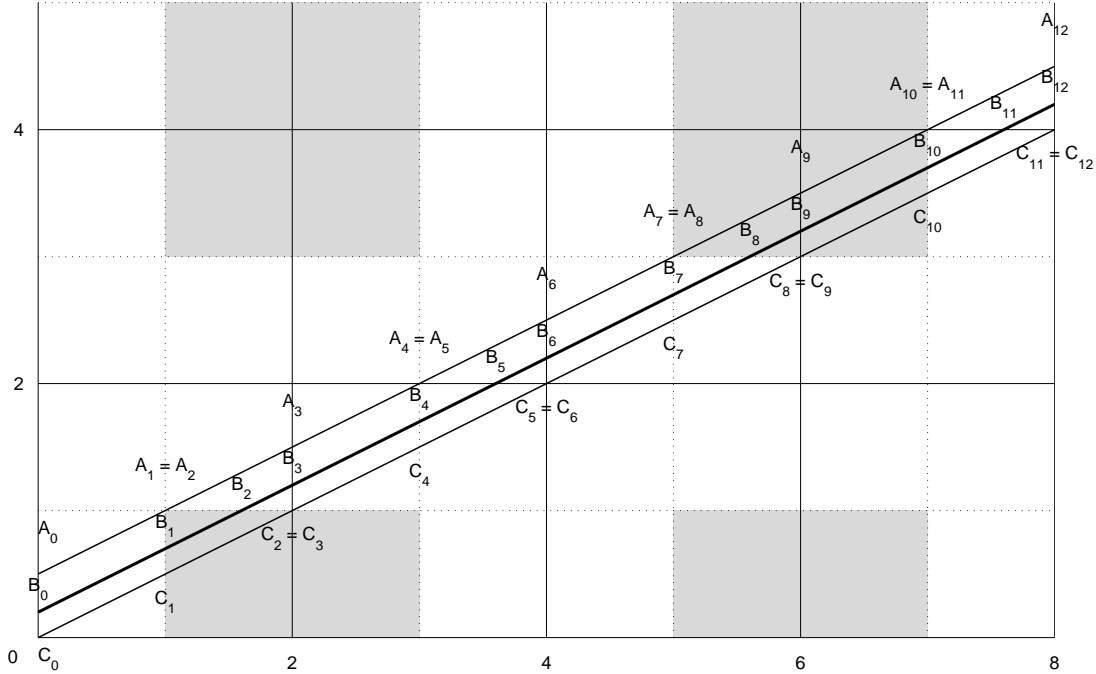,height=9cm}} 
\caption{Case of a ray which does not meet any point with integer coordinates} \label{fig_points_nonentiers} 
\end{figure} 
 
\noindent Let $A_i$ be the intersection point of the line meeting
$A$ and 
 directed by 
$\vec{u}$ (the points  $A$ solutions of the minimization problem above are
clearly on a same line directed by $\vec{u}$) and  the line belonging to
the lattice $\Z+\Z$ which contains  $B_i$ (see 
Figure \ref{fig_points_nonentiers} for the case  $N=2$). In the same way, 
there exists  $C$  with integer coordinates which satisfies 
$$d(C,\gamma)= 
\min \{ d(P,\gamma) | P \text{ has integer coordinates and } d(P,\gamma)<0
\}. $$

\noindent Let us note $C_i$ 
the intersection  of the line meeting  $C$ and
directed by 
$\vec{u}$ and the line belonging to
the lattice $\Z+\Z$ which contains  $B_i$  (see Figure 
\ref{fig_points_nonentiers}).\\ There are  three oblique rays of angle 
$\alpha$ $\gamma_A = 
[A_0,\alpha]$, $\gamma_B = [B,\alpha] = [B_0,\alpha]$ and 
$\gamma_C = [C_0,\alpha]$. According to  the definitions  of $A$ and $C$, the 
segments $[A_i\,A_{i+1}]$, $[B_i\,B_{i+1}]$ 
and  $[C_i\,C_{i+1}]$ belong to the same square
$C_{k,l}$. We define  $\varepsilon_k=1$  if
$[B_k\,B_{k+1}] \subset \om_T$, $\varepsilon_k=0$ if $[B_k\,B_{k+1}] \not\subset \om_T$. Then, if 
$$m(\gamma_B)=\frac{1}{L} \int_{[B_0 
B_{2N(p+q)}]} \chi_{\om_T} =  \frac{1}{L} \sum_{k=1}^{2N(p+q)} 
\varepsilon_k \, B_{k-1} B_k,$$ 
we have: 
$$m(\gamma_A)=\frac{1}{L} \sum_{k=1}^{2N(p+q)} \varepsilon_k \, 
A_{k-1} A_k \quad \text{and} \quad m(\gamma_C)=\frac{1}{L} 
\sum_{k=1}^{2N(p+q)} \varepsilon_k \, C_{k-1} C_k.$$ These three rays are
parallel. Hence, there exists $t\in(0,1)$ such that for all  
 $i$, $B_i$ is the  barycenter of  $(A_i,t)$ and  $(C_i,1-t)$. If the
 quadrilateral  
$A_i\,A_{i+1}\,C_i\,C_{i+1}$ is a trapezium (or in the degenerate case a
rectangular triangle), 
Thales' theorem implies that $B_i B_{i+1} = t \, A_i 
A_{i+1} + (1-t)\, C_i C_{i+1}$. If the quadrilateral is a
parallelogram, this inequality remains valid because the three lengths
are equal. As a consequence, we have $$ m(\gamma_B) = t 
m(\gamma_A) + (1-t) m(\gamma_C).$$ The minimum of the three quantities
$m(\gamma_A)$, $m(\gamma_B)$ and $m(\gamma_C)$ is attained for one of the
two rays 
 $\gamma_A$ ou $\gamma_C$. This shows that $\gamma$ cannot be
 the minimum in the set of oblique closed rays. \fpr\\ 
 
\noindent It remains to find the angles $\alpha$ we need to consider  among the rays which meet a point  with
integer coordinates. In the flat torus
$T=[0,\, 2N]^2$, there are 
$4N^2$ possible starting points.  However, the ray  $\gamma 
=[(x_0,y_0),\alpha]$ with  $(x_0,y_0) \in \N^2$ meets $2N$ 
points with integer coordinates. Indeed, $\gamma$ 
meets the point $(x_1,y_1)$ if and only if there exists $\l 
\in [0,2N)$ such that $$\left\{ 
\begin{array}{ll} x_1=x_0+\l p\; \;[2] \\y_1=y_0+\l q \; \; [2]\end{array} 
\right. .$$ Therefore, the points belonging to $\gamma$ and 
 with integer coordinates are exactly the  $2N$ points obtained for
$\l = 
0,1,\dots, 2N-1$. There exist  at least $2N$ rays which meet a point
 with integer coordinates. 
\begin{proposition} \label{th_pointsdep} 
Let $d$ be the G.C.D.  of $p$ and $2N$, let 
$d'$ be such that $d\,d'=2N$. Then, the  $2N$ rays  of angle $\alpha$ 
starting at $(i,j)$ with $0 \leq 
i \leq d-1$ and $0 \leq j \leq d'-1$ meet one and only one time all the
points of $T$  with integer coordinates.
\end{proposition} 
\noindent \Pr Let  $p'$ be  such that $d\,p'=p$. We note that  $p'$ and 
$d'$ are relatively prime. Let also $\gamma_1$ and 
$\gamma_2$ be the rays starting respectively at $A$ and $B$ 
with integer coordinates $(k_1,k_2)$ and $(k_3,k_4)$ (with  $0 
\leq k_1,\, k_3 \leq d-1$ and $0 \leq k_2,\, k_4 \leq d'-1$ ).
There exists a point $C \in \gamma_1 \cap \gamma_2$  with integer
coordinates
 if and only if there exist four integer numbers
$\l,\mu,n$ and $m$ such that $$\left\{ 
\begin{array}{ll} x = k_1 + \l p = k_3 + \mu p+ 2 N n\\ 
  y = k_2 + \l q = k_4 + \mu q + 2 N m\end{array} \right. .$$ Therefore, 
$k_1-k_3 = (\mu - \l)p + 2Nn$ and $k_1 - k_3$ is divisible by $d$. 
Since $-(d-1) \leq k_1 - k_3 \leq d-1$, we have $k_1=k_3$ and $(\mu - \l)p 
= -2N n$. As a consequence, $(\mu - \l)p' = -d' n$ and since $p'$ and 
$d'$ are relatively prime, $(\mu - \l) = k d'$, where $k \in \Z$. In the
same way, $k_2-k_4 = (\mu - \l)q + 2N m = kd'q + 2N m$. 
Hence $k_2 - k_4$ is divisible by $d'$. Since  $-(d'-1) \leq k_2 - k_4 
\leq d'-1$, we have  $k_2=k_4$ and $\gamma_1 = \gamma_2$. \fpr

\subsection{The algorithm} 
In this section,  we give a method  to compute
explicitly   $g$ when $\om$ is a finite union of squares $C_{i,j}$. 
At first, according to Propositions 
\ref{prop1} and  \ref{th_points_nonentiers}, 
the study can be restricted to closed rays 
starting at a point  with integer coordinates. \\ The user of the program
must input the matrix $M$ which 
represents $\om$ and the value of a parameter  $PQmax$ which 
avoids infinite loops (this problem   never appeared
until now).\\ 
 
\noindent The algorithm is the following 
\begin{enumerate} 
\item  Compute  the minimum $g$ of the $2N$ numbers $$\frac{1}{N} \sum_{i=1}^N M_{i,j} 
\; \text{ and } \; \frac{1}{N} \sum_{j=1}^N M_{i,j},$$ which corresponds to  
the  minimum of $m(\gamma)$ when $\gamma$ is a vertical or horizontal ray.

\item From the matrix $M$, build the matrix $A$ defined in  
(\ref{eq_defA}). 
\item For all $(p,q)$ such that  $p \geq 1$, $q \geq 1$, $p+q \leq 2N$ and $p$ 
and $q$ relatively prime, compute 
the number (the way to compute $m_{p,q}$ is explained below) 
\begin{equation} \label{eq_def_mpq} 
 m_{p,q} = \min_{X \in \{1,2,\dots,2N\}^2} \{ m(\gamma) | \gamma = [X,
 \arctan q/p] \}.
\end{equation} 
Then $g \leftarrow  \min \{g, m_{p,q}\}.$ 
\item If at this step $g=|\om|$, the user inputs the parameter $PQmax$ 
which corresponds to the greatest product $p q$ which will be considered~; 
in the
other cases $PQmax = E(\frac{N^2}{|\om| - g} \min\{|\om|,1-|\om|\})+1$. 
\item Find all the  $(p,q)$ relatively prime such that $p \geq 1$, $q \geq 1$, $p+q > 2N$, 
$pq \leq PQmax$. Sort these
couples from the lowest product $pq$ to the greatest  one. 
This gives a list $L$ containing $n$ couples. 
\item $i \leftarrow 1$ and as long as $|\om| - g \geq  \frac{N^2}{pq} \min \{ |\om|,1-|\om| \}$ and 
$i \leq n$, $$ g \leftarrow  \min \{g, m_{p,q}\} \; \text{ and } \; 
i \leftarrow i+1,$$ where $(p,q)$ is the  {\it i-th} couple of 
 $L$ and $m_{p,q}$ is defined above (see (\ref{eq_def_mpq})). 
\item Finally 
\begin{itemize} 
\item if at the end of this loop, $g = |\om|$, then$$|\om| - \frac{N^2}{PQmax} \min\{|\om|,1-|\om|\} \leq 
g(\om) \leq |\om|.$$ This means that one of the following cases occurs;
$g(\om) = |\om|$, or too few  families of rays have been considered.
\item if $g \not= |\om|$, according to  Theorem \ref{th_influencealpha}, $g(\om)$ is
  equal to $g$.  
\end{itemize} 
\end{enumerate} 
 
\vspace{0.5cm} 
 
\begin{figure}[h] 
\centerline{\epsfig{figure=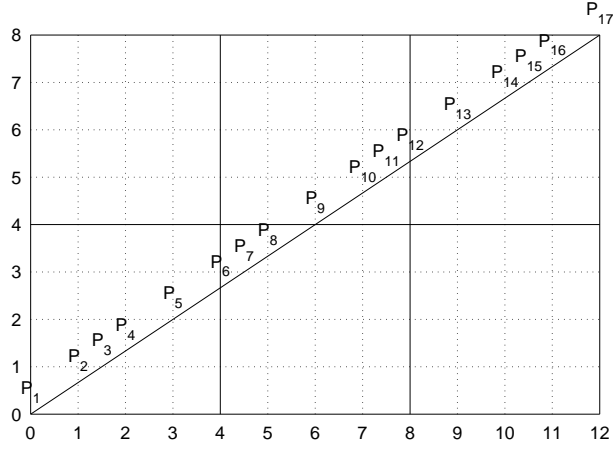,height=6cm}} 
\caption{Ray $\gamma_0$, for $N=2$, $p=3$ and $q=2$} 
\label{fig_mpq} 
\end{figure}

We now present the method of computation defined in
(\ref{eq_def_mpq}) for two relatively prime integers 
 $p \geq 1$, $q \geq 1$. At first,  
let us consider the ray  $\gamma_0 = [0, 
\arctan q/p]$ and let us study the way to compute $m(\gamma_0)$. Its period
is $L= 
2N\sqrt{p^2+q^2}$ and between the instants $t=0$ and  $t=L$, it meets  
$2N(p+q-1)+1$ points which have at least one integer coordinate.
Let $P_1$,..., $P_{2N(p+q-1)+1}$ be these points (see 
Figure \ref{fig_mpq}). Remember that  $P_{l(p+q-1)+1}$, 
$1 \leq l \leq 2N$ have integer coordinates. \\ 
 
\noindent Let $(i_k)$, $(j_k)$ and $(\varepsilon_k)$ be the three other finite sequences
($1 
\leq k \leq 2N(p+q-1)$) defined by $$P_k P_{k+1} \subset 
[i_k,i_k+1] \times [j_k,j_k+1]$$ and $$ \varepsilon_k = \left\{ 
\begin{array}{ll} 1 \text{ if } [i_k,i_k+1] \times [j_k,j_k+1] 
\subset \om\\ 0 \text{ if not } \end{array} \right. .$$ Then, 
$m(\gamma_0) = \frac{1}{L} \sum_{k=1}^{2N(p+q-1)} \varepsilon_k 
P_k P_{k+1}$, the computation of  $m(\gamma_0)$ is equivalent to the
computation of the sequences  $(P_k P_{k+1})$, $(i_k)$ and $(j_k)$.\\ 
 
\noindent $\bullet$ Computation of the lengths $P_k P_{k+1}$\\ At first,
note that it suffices to compute the lengths
$P_k P_{k+1}$ for $1 \leq k \leq p+q-1$. Indeed, 
$\overrightarrow{P_{k+(p+q-1)}P_k} = p \vec{i} + q \vec{j}$, hence 
the sequence $(P_k P_{k+1})$ is periodic of period $p+q-1$. Since 
 $p$ and $q$ are relatively prime,  $P_1$  and 
$P_{p+q-1}$ do not have integer coordinates.\\ 
 
\noindent Then, the equations of the ray 
$\gamma_0$ are $x = p t/ \sqrt{p^2+q^2}$ and $y = q t/ 
\sqrt{p^2+q^2}$. It meets  the vertical lines defined by 
$x = i$ at  instants $t_i =  i \sqrt{p^2+q^2}/p $ 
at points $(i, i q/p)$ and the horizontal lines defined by 
$y=j$ at instants $t'_j =  j \sqrt{p^2+q^2}/q $ and at points $(j 
p/q, j)$.\\ Let us then define the list $V^1$ by $$V^1 = 
\{t'_j,1 \leq j \leq q-1\} \cup \{t_i,1 \leq i \leq p-1\}.$$ Let 
$V^2$ be the list obtained from $V^1$ by sorting its elements in
increasing order. Let also $\sigma$ be the permutation  of $\{1,2, \dots, 
p+q-2\}$ such that $$V^2 = \sigma(V^1) = \{V^1_{\sigma(1)}, 
V^1_{\sigma(2)}, \dots, V^1_{\sigma(p+q-2)}\}.$$ Let $V^3$ and 
$V^4$ be defined by $$V^3 = \{0\} \cup V^2 \cup 
\{\sqrt{p^2+q^2}\} \; \text{ and } \; V^4 = \{V^3_2 - V^3_1, V^3_3 
- V^3_2, \dots, V^3_{p+q} - V^3_{p+q-1}\}.$$ Finally, let $V$ be defined by 
$$\begin{array}{cc} V = & \underbrace{V^4 \cup 
V^4 \cup \dots \cup V^4}\\ & 2N \text{ times} \end{array},$$ then for all  
$k$, $1 \leq k \leq 2N(p+q-1)$ $P_k P_{k+1} = V_k$ 
and $$m(\gamma_0) = \frac{1}{2N\sqrt{p^2+q^2}} 
\sum_{k=1}^{2N(p+q-1)} \varepsilon_k V_k.$$ 
 
\noindent As an important remark, it may be convenient to work with  $pq /
\sqrt{p^2+q^2}\, .V$ instead of  
$V$. Thus,  
$V_1 = \{p j,1 \leq j \leq q-1\} \cup \{q i,1 \leq i 
\leq p-1\}$, $V_3 = \{0\} \cup V^2 \cup \{pq\}$ and hence, the elements of
$V$ are all integer numbers. $m(\gamma_0)$ is then given by 
$$ m(\gamma_0) = \frac{\sum_{k=1}^{2N(p+q-1)} 
\varepsilon_k V_k}{2N pq}.$$ As a consequence, $m(\gamma_0)$ 
is a rational number whose denominator and numerator are explicitly known.\\ 
 
\noindent $\bullet$ Computation of the sequences $(i_k)$ and $(j_k)$\\ Let 
$I^1$ and $J^1$ be the two following lists $$\begin{array}{cccccc} 
I^1 = & \{ \underbrace{1,1,\dots,1}, & \underbrace{ 0,0,\dots,0} 
\} & \; \text{ and } \; J^1 = & \{ \underbrace{0,0,\dots,0}, & 
\underbrace{1,1,\dots,1}\} \\ & q-1 \text{ times} & p-1 \text{ 
times}& & q-1 \text{ times} & p-1 \text{ times} \end{array}.$$ 
This means that the numbers 1 which appear in the sequence $I^1$ correspond
to the points 
$P_k$ which belong  to horizontal lines and the numbers 1 of the sequence 
$J^1$ correspond to the points $P_k$ belonging to vertical lines. 
From these two lists and from the permutation $\sigma$, we define 
$$ I^2 = \sigma(I^1) = \{I^1_{\sigma(1)}, 
I^1_{\sigma(2)},\dots, I^1_{\sigma(p+q-2)}\} \; \text{ and } \; J^2 
= \sigma(J^1) = \{J^1_{\sigma(1)}, J^1_{\sigma(2)},\dots, 
J^1_{\sigma(p+q-2)}\}$$ and also $$ I^3 = \left\{0, I^2_1 , 
\sum_{k=1}^2 I^2_k, \dots, \sum_{k=1}^{p+q-2} I^2_k\right\} \; 
\text{ and } \; J^3 = \left\{0, J^2_1 , \sum_{k=1}^2 J^2_k, \dots, 
\sum_{k=1}^{p+q-2} J^2_k\right\}.$$ Finally, we set $$I = I^3 \cup 
(p+I^3) \cup \dots \cup ([2N-1]p + I^3) \; \text{ and } J = J^3 
\cup (q+J^3) \cup \dots \cup ([2N-1]q) + J^3.$$ Then, for all  
$k$, $1 \leq k \leq 2N(p+q-1)$, $i_k = 1+I_k$ and $j_k = 1+J_k$.\\ 
Let us  introduce the matrix $A$ defined from $M$ 
in (\ref{eq_defA}) and a function $r$ defined over integers by 
$$r(k) \in \{1,2, \dots, 2N\}\; \text{ and } \; \exists \, l 
\in \Z \,|\, k= 2N.l + r(k).$$ Then, for all $k$, $1 \leq k \leq 
2N(p+q-1)$, $\varepsilon_k = A_{r(i_k),r(j_k)} 
=A_{r(1+I_k),r(1+J_k)}$.\\ 
 
\noindent There are 
 two reasons why we consider all rays of angle $\arctan q/p$ 
at the same time in the computation of $m_{p,q}$. First, all these rays
have the same sequence  $V_k$. In addition, let $\gamma= [X_0,
\arctan (\frac{q}{p}) ]$ where $X_0=(x_0,y_0)$ has integer coordinates. 
Then, the sequences $(i_k)$ and 
$(j_k)$ of $\gamma$ can be computed from the lists  $I_k$ 
and $J_k$ by the formula: $i_k = 
1+x_0+I_k$ and $j_k = 1+x_0+J_k$.\\ 
 
\noindent For two relatively prime integers $p \geq 1$ and $q \geq 1$, the 
number $m_{p,q}$ can be computed in the following way 
\begin{enumerate} \item compute the G.C.D. 
$d$ of $p$ and  $2N$ and from Proposition  
\ref{th_pointsdep}, the  $2N$ starting points  with integer
coordinates  that must be considered are known. Let us  note $X_1, X_2, 
\dots, X_{2N}$ and  $Y_1, Y_2, \dots, Y_{2N}$ their coordinates. 
\item Compute the three lists of  $2N(p+q-1)$ elements 
$V_k$, $I_k$ and $J_k$ defined above. \item $m_{p,q}$ is then the minimum
of the $2N$ numbers $$ 
\frac{\sum_{l=1}^{2N(p+q-1)} A_{r(1+X_k+I_l),r(1+Y_k+J_l)}\, 
V_l}{2N pq}.$$ All these quotients have the same denominator and hence,
are easy to compare by looking at their numerator. Finally, this gives
the explicit value of $m_{p,q}$. 
\end{enumerate} 
 
\vspace{0.5cm} 
 
\noindent We now present the results obtained  with our algorithm for explicit domains
$\om$.

\begin{figure}[h] 
\centerline{\psfig{figure=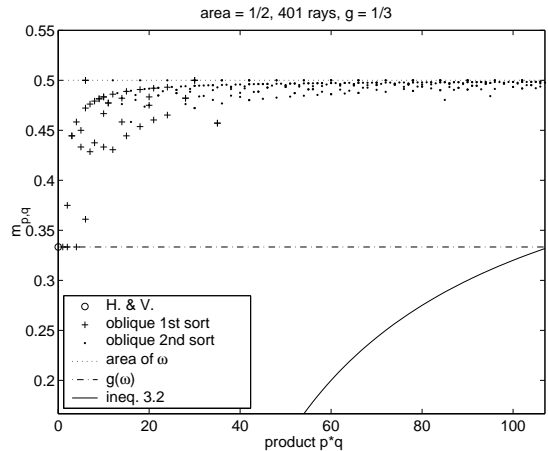,height=6cm}} 
\caption{An example of computation of $g(\om)$} 
\label{fig_exemple_calcul} 
\end{figure}

\begin{figure}[p] 
\centering 
 a.\epsfig{figure=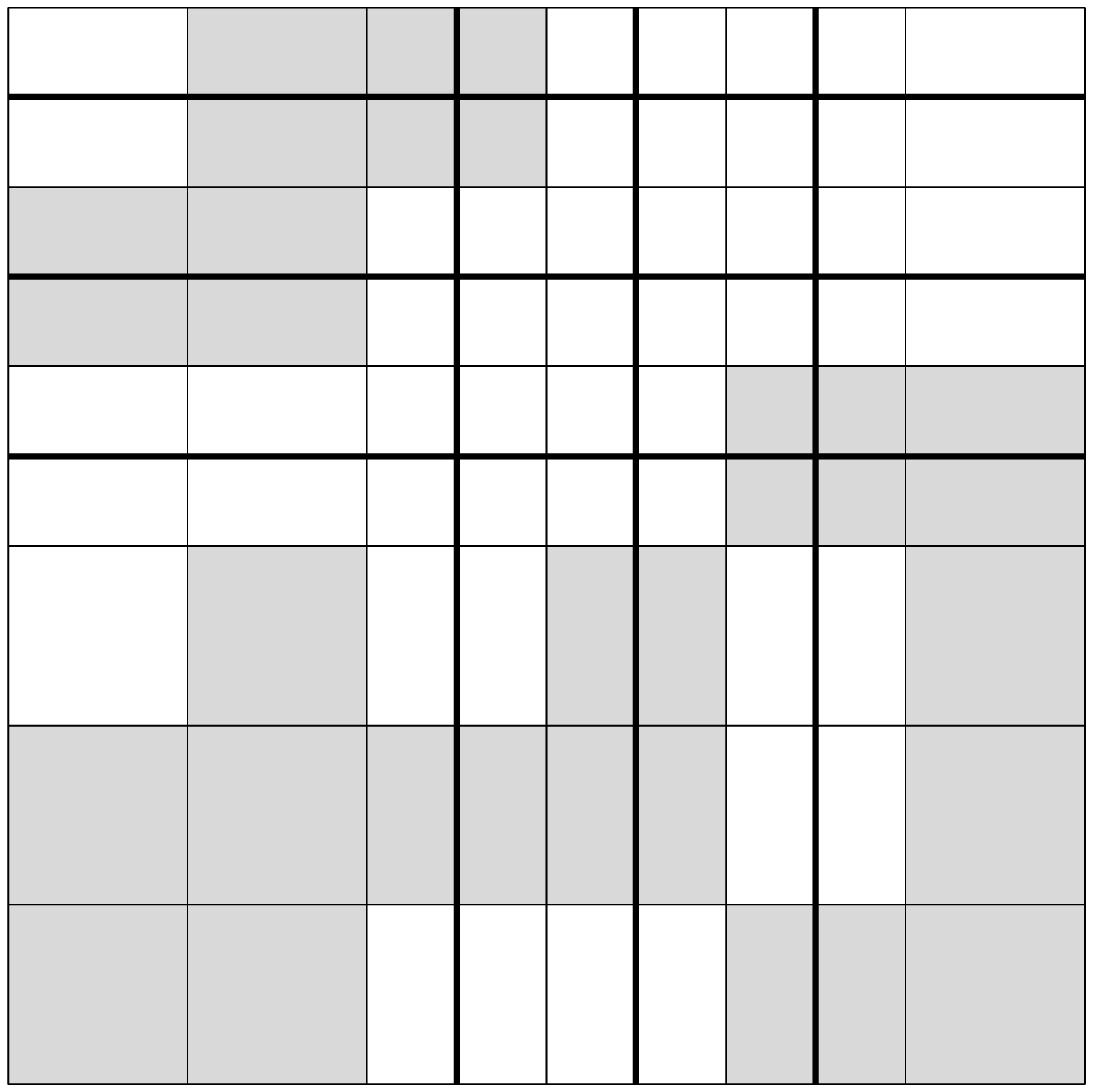,height=4cm} 
~b.\epsfig{figure=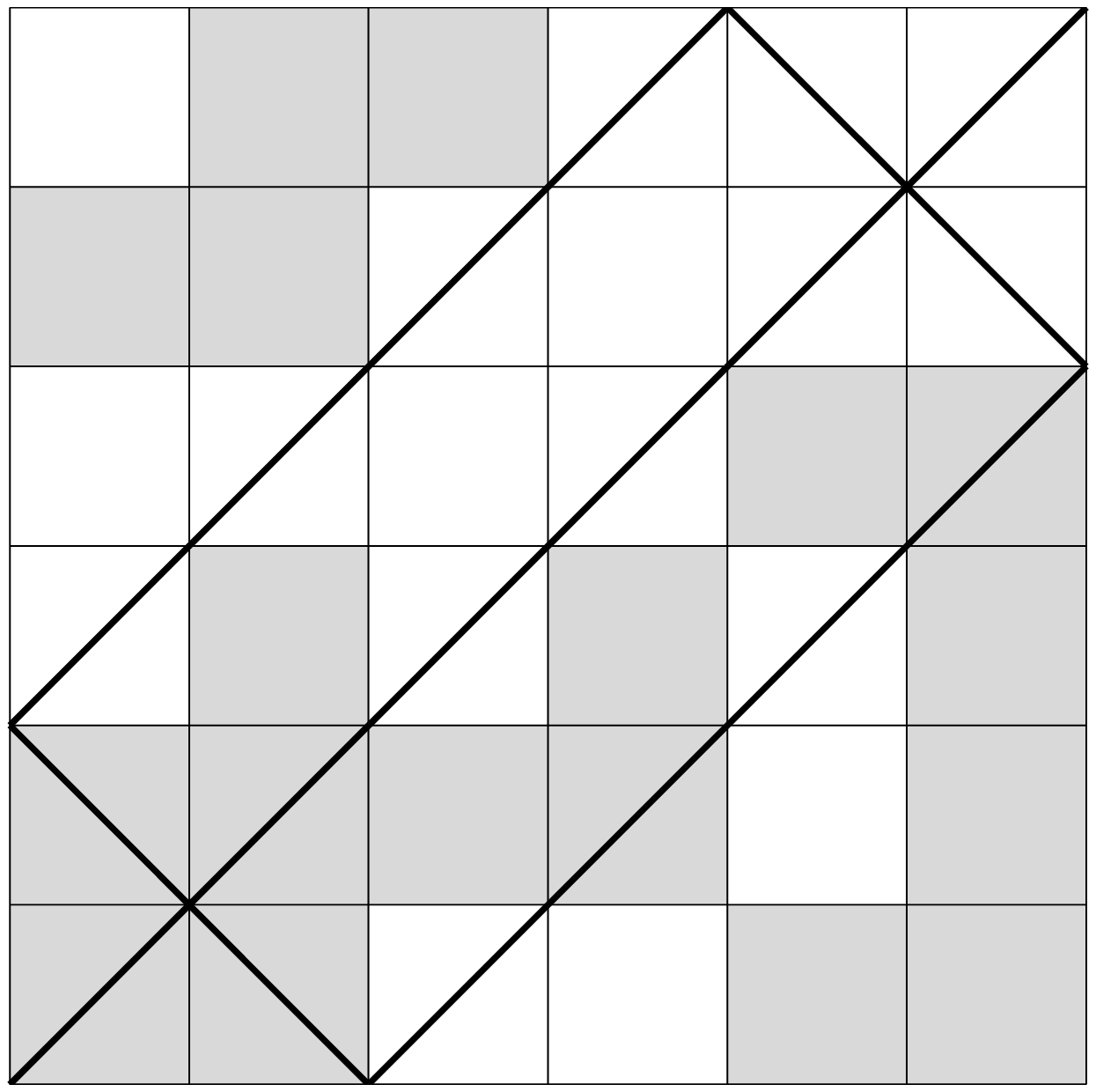,height=4cm}\\ 
 c.\epsfig{figure=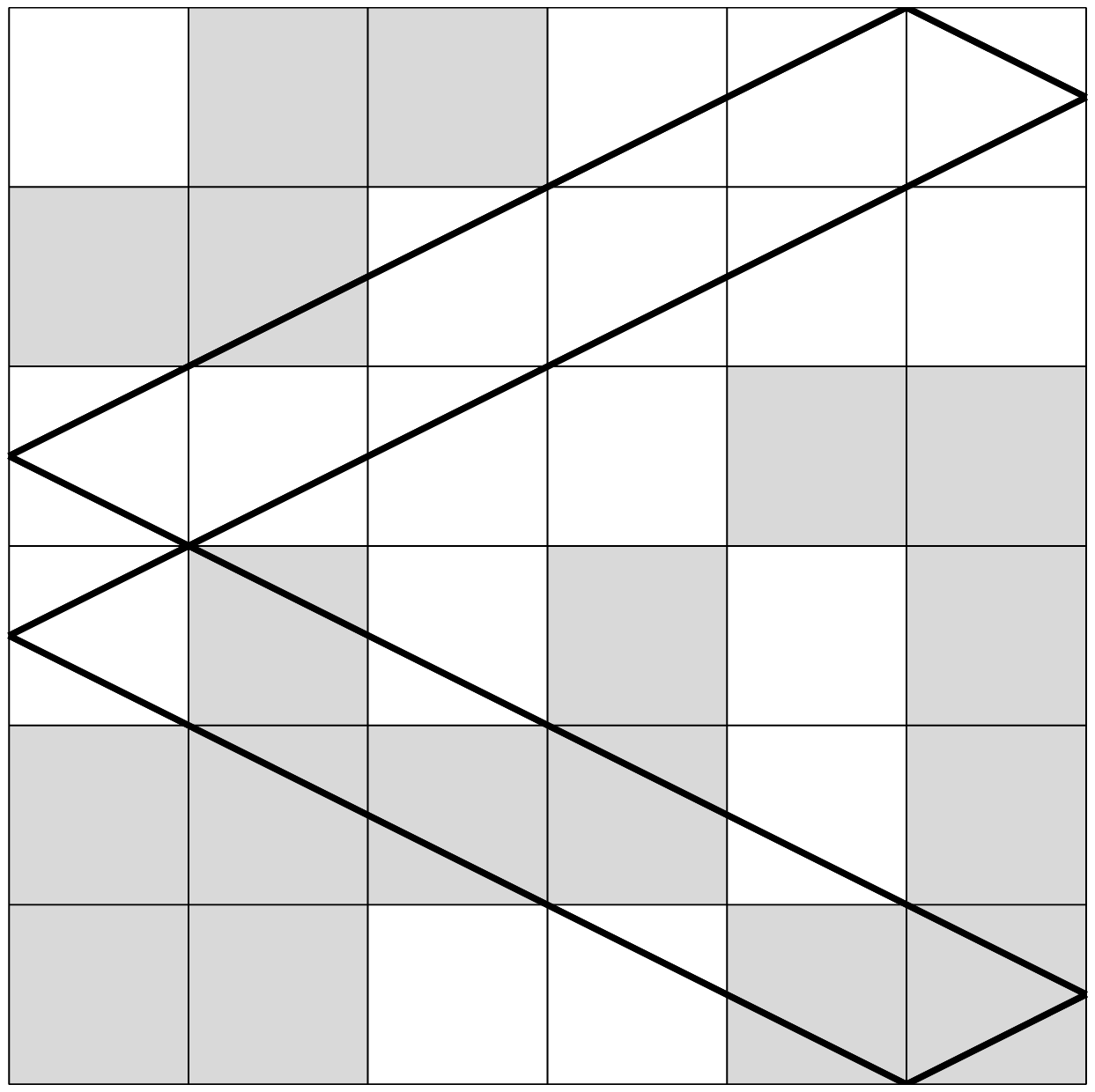,height=4cm} 
~d.\epsfig{figure=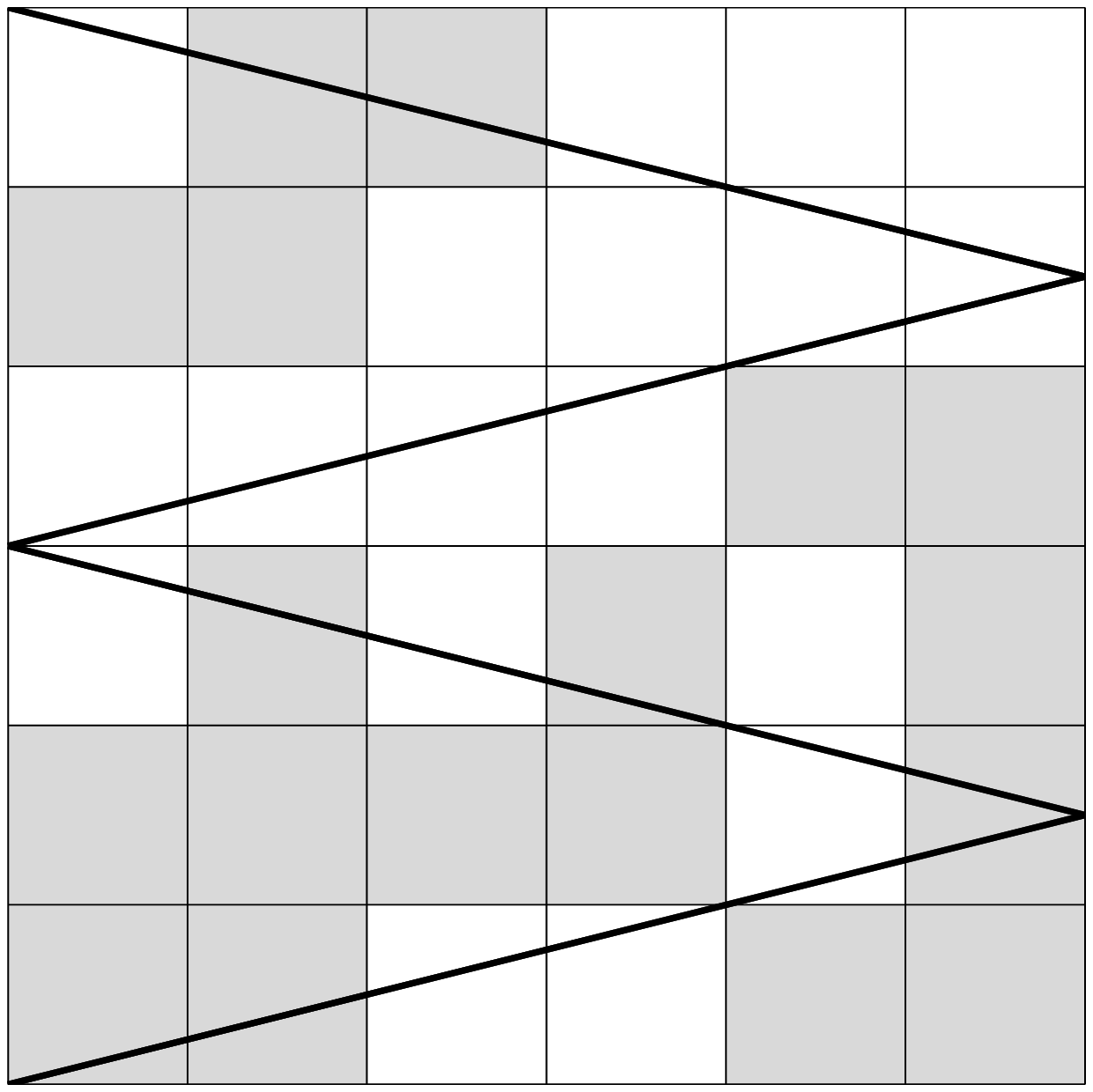,height=4cm} 
 
\caption{Minimizing rays : $m(\gamma) = g(\om) = 1/3$} 
\label{fig_raysat1} 
\end{figure} 
 
\begin{figure}[p] 
\centerline{\epsfig{figure=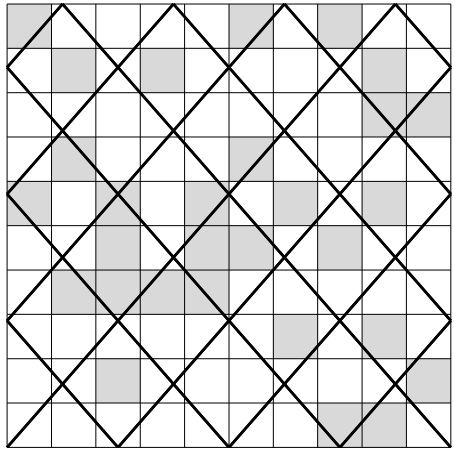,height=6cm}} 
\caption{Minimizing ray : $m(\gamma) = g(\om) = 111/560$} 
\label{fig_raysat2} 
\end{figure} 
As an example, for the computation of $g(\om)$ for the domain $\om$ of Figure
\ref{fig_codage}, the minimum for 
vertical and horizontal rays is  $g=1/3$ and 
there are  45 families of oblique rays of first kind 
(i.e. the rays of part 3 in algorithm). Among these rays,   
there exists  $\gamma$ such that $m(\gamma)=1/3$ but there is no $\gamma$
such that  
$m(\gamma)<1/3$. At this step, the value of  $g$ is $1/3$. 
There are also  344 families of oblique rays of second kind (corresponding to
part 5 of algorithm). None of them satisfies
$m(\gamma)=1/3$. We obtain $g(\om) = 1/3$.\\

\noindent Figure \ref{fig_exemple_calcul} presents the values of 
$m_{p,q}$ as a function of $pq$ and presents also inequality of Theorem 
\ref{th_influencealpha}. Note that this inequality is not  
optimal for this particular domain. However, (\ref{eq_encadcarre}) is optimal for a
unique square 
$C_{i,j}$. This shows that the error comes from the summation of
inequalities (\ref{eq_encadcarre}) for all the squares of
$\om_T$. Nevertheless, it  is useful because it certifies  that the
value of $g(\om)$ is exact.\\
 
\noindent There exist 10 rays which satisfy  $m(\gamma) = g(\om)$; six 
vertical or horizontal rays (more precisely  six  families of vertical or
horizontal rays). They are represented on Figure 
\ref{fig_raysat1}.a. There are also  two rays of angle $\pi/4$  
represented on Figure \ref{fig_raysat1}.b, one ray
of angle $\arctan 1/2$ represented on Figure 
\ref{fig_raysat1}.c, and also the  ray $\gamma = [0, \arctan 
1/4]$ represented on Figure \ref{fig_raysat1}.d.\\ 
 
\noindent The domain $\om$ of Figure \ref{fig_raysat2} of area  
$|\om| = 29/100$ is represented here to show that intuition can be
false. Indeed, consider only vertical rays, horizontal rays and
closed rays of small period. The intuition says that 
 $g(\om) = 1/5$. However, the ray  $\gamma = [0, \arctan 
8/7]$ is such that $m(\gamma) = 111/560$ (note that  1/5=112/560), 
that is the exact value of $g(\om)$ as shown by the complete computation.

\section{Appendix : proof of Theorem \ref{th_egcrit}} 

\noindent  First, it is clear that 
\begin{eqnarray} \label{r1}
g(\om) \leq g'(\om).
\end{eqnarray}
Indeed, for all $\gamma' \in \mathcal{C}$ and all $t >0$,  
$m(\gamma_t' ) \geq  \inf_{\gamma \in \mathcal{C}} m(\gamma_t)$. Hence, 
$$m(\gamma')= \limsup_{t \to +\infty} m(\gamma'_t) \geq 
\limsup_{t \to +\infty} \inf_{\gamma \in \mathcal{C}} m(\gamma_t)=
g(\om).$$
Since this equality is 
true for all  $\gamma' \in \mathcal{C}$, this shows 
(\ref{r1}).

\noindent Let us  prove that 
\begin{eqnarray} \label{r2}
g(\om) \geq g'(\om).
\end{eqnarray}
Let us fix $\epsilon >0$ and 
choose a sequence 
$(\gamma^n)_n \subset \mathcal{C}$ and a sequence $(t_n)_n$  which tends to
$+\infty$ such that

\begin{eqnarray} \label{r3}
\lim_{n \to +\infty} m(\gamma^n_{t_n})= 
g(\om).
\end{eqnarray}

\noindent For all $n \in \N$, there exist $X_n \in \Om$ and $\theta_n \in 
[0,2 \pi[$ such that $\gamma^n = [X_n,\theta_n]$. Up to a 
subsequence, we can assume that there exists $(X_0,\theta) 
\in \overline{\Om} \times [0,2 \pi]$ such that 
$$\lim_n X_n=X_0  \hbox{ and }  \lim_n \theta_n= \theta.$$
We claim that there exists $L>0$ such that 
\begin{eqnarray} \label{r5}
m([y,\theta]_L) \geq g'(\om)-\epsilon  \hbox{ for all } y \in \Om. 
\end{eqnarray}
If $\tan(\theta) \in \R\backslash\Q$, then set  $L=t_0$, where $t_0$
is as in Proposition  \ref{prop1}. 
Indeed, together with relation  (\ref{b1}), 
$$ m([y,\theta]_L) \geq |\om|-\epsilon
\geq g'(\om)-\epsilon.$$
If $\tan(\theta) \in \Q$ or if $\tan(\theta)=\infty$,
the ray $[y,\theta]$ is periodic and the period $T$ 
does not depend on $y \in
\Om$. 
Hence,
$$m([y,\theta]_L)=m([y,\theta]) \geq g'(\om).$$
Thus, let us set $L=T$. 
Clearly, $t_n$ can be assumed to be a multiple of $L$. Indeed, replacing 
$t_n$ by $t_n+s_n$ (where $|s_n| \leq L$) does not change the limit of 
$m(\gamma^n_{t_n})$. Now, assume that  $\theta_n = \theta$.
If $\tan(\alpha) \in \Q$ or if $\tan(\alpha)=\infty$,
then clearly $\lim_n m(\gamma^{n}_{t_n})- m(\gamma^n)=0$. If 
$\tan(\alpha) \in \R \backslash \Q$, this relation follows from Proposition
\ref{prop1}. Then, by relation (\ref{r3}), it follows  that 
$$g(\om) \geq \liminf_n  m(\gamma^n) \geq g'(\om).$$
This proves (\ref{r2}). Hence, we can assume that 

\begin{eqnarray} \label{r1'}
\theta_n \not= \theta.
\end{eqnarray}

\noindent Now, let us define, for $k \in \{ 
0,...,\frac{t_n}{L}-1 \}$,
$X_{k,n}=\gamma^n(kL)$ and $\gamma^{k,n}(t)= \gamma^n(kL+t)$ for 
$t\in [0,L]$. 
Let us define also 
${\gamma}^{k,n,s}$ as the unique trajectory of length $L$ and slope
$\theta$  such that 
${\gamma}^{k,n,s}(s)= \gamma^{k,n}(s)$
and 
$$A_{k,n}= \bigcup_{s \in [0,L]} {\gamma}^{k,n,s}.$$


 \begin{figure}[h] 
\centerline{\input{carr1.pstex_t}}
\end{figure}



\noindent Now, we prove that

\begin{step} 
There exists a sequence $(k_n)_n$ such that 
$$|A_{k_n,n} \cap \partial \omega| \leq \epsilon \hbox{ and } 
m(\gamma^{k_n,n}) \leq g(w)+\epsilon $$
for all $n$ large enough. Here, $|A_{k_n,n} \cap \partial \omega| $ denotes
the length of the curve $A_{k_n,n} \cap \partial \omega$.
\end{step}

\noindent In the following, the length of any curve $C$ will be denoted by
$|C|$. Let 
$$\mathcal{K}_n= \left\{ k \in   \{ 
0,...,\frac{t_n}{L}-1 \} \hbox{ s.t. } m(\gamma^{k,n}) \leq g(w)+\epsilon \right\}.$$
Set
$$\mathcal{M} = \bigcap_{m \in \L} \left(\overline{\bigcup_{n \geq m; 
k \in \mathcal{K}_n } X_{k,n}}
\right)$$  
and 
$$\tilde{\mathcal{M}} = \bigcap_{m \in \L} \left(\overline{\bigcup_{n \geq m ;
k \in \mathcal{K}_n } \gamma^{k,n}}
\right).$$  
For $X \in \mathcal{M}$, let $\gamma_X$ denote the trajectory 
of length $L$ and slope $\theta$ such 
that $\gamma_X(0)=X$.
Note that $\mathcal{M}$ and $\tilde{\mathcal{M}}$ are not empty. Otherwise,
$\mathcal{K}_n$ would be empty too and we would 
have $m(\gamma^{k,n}) \geq  g(\om)+ \epsilon$ for all $n \geq n_0$ and forall
$k$. Thus, we would have $  m(\gamma^n_{t_n})\geq  g(\om)+ \epsilon$ that 
contradicts (\ref{r3}). 
Now, we prove that

\begin{eqnarray} \label{r6}
\tilde{\mathcal{M}} = \bigcup_{y \in \mathcal{M}} \gamma_y.
\end{eqnarray}
Take a sequence $(X_{l_n,n})_n$ such that
$\lim_n X_{l_n,n}=y \in \mathcal{M}$. Then, 
for all $t \in [0,L]$,
$\lim \gamma^{l_n,n}(t)= \gamma_y (t)$. We obtain that
$$\bigcup_{y \in \mathcal{M}} \gamma_y  \subset \tilde{\mathcal{M}}. $$
Conversely, if $z \in \tilde{\mathcal{M}}$ then there exist two  sequences 
$(l_n)_n$ such that $l_n \in \mathcal{K}_n$ and $(s_n)_n \subset
[0,L]$ such that $z=\lim_n \gamma^{l_n,n}(s_n)$. Up to a subsequence, we can 
assume that there exist $y \in \mathcal{M}$ and $s \in [0,L]$ such that 
$\lim_n X_{k_n,n}=y$ and $\lim_n s_n=s$.  As one can check, this implies that
$z= \lim_n \gamma^{k_n,n}(s)= \gamma_y(s)$ and therefore 
$$ \tilde{\mathcal{M}}\subset \bigcup_{y \in \mathcal{M}} \gamma_y.$$
This proves (\ref{r6}).

\noindent We now distinguish two cases.\\

\noindent First, assume that
there exists an infinite sequence $(y_n)_n$ in $\mathcal{M}$ such that 
$\gamma_{y_n} \cap \gamma_{y_{n'}} = \emptyset$ if $n \not= n'$.  
Since $|\partial \omega|$ is finite, one can choose $y \in \mathcal{M}$ ($y$ 
is one of the $y_n$) 
such that $|\gamma_y \cap \partial \omega | \leq \frac{\epsilon}{2} $.
Otherwise, for all $n$, $|\partial \omega \cap \gamma_{y_n} | \geq
\frac{\epsilon}{2}$ and hence, since  $\gamma_{y_n} \cap \gamma_{y_{n'}} =
\emptyset$, we would have 
$$+ \infty = \sum_{n} |\partial \omega \cap \gamma_{y_n} | = 
\left| \left( \cup_n \gamma_{y_n} \right) \cap \partial \omega \right| \leq
|\partial \omega|.$$
This is false.
  The definition of $\mathcal{M}$ implies that there exists a sequence 
$(k_n)_n$ such that $y= \lim_n X_{k_n,n}$ with $k_n \in \mathcal{K}_n$.
Clearly, for all $t \in [O,L]$,

\begin{eqnarray} \label{r7}
\lim_n \gamma^{k_n,n}(t) = 
\gamma_y(t).
\end{eqnarray}
Set $B_n= \overline{\bigcup_{m \geq n} A_{k_n,n}}$. 
Obviously, $B_n$ is a decreasing sequence of sets. We recall that 
$$A_{k_n,n}= \bigcup_{s \in [0,L]} \gamma^{k_n,n,s}$$
By construction, for all $t \in [0,L]$, $ \gamma^{k_n,n,s}(t)=
\gamma_y(t)$.
This implies that $\bigcap_{n \in \N} B_n= \gamma_y$. Since 
$|\gamma_y \cap \partial \omega | \leq \frac{\epsilon}{2} $, for all 
$n$ large enough, we have $|B_n \cap \partial \omega| \leq \epsilon$. Since
$A_{k_n,n} \subset B_n$ and since $k_n \in \mathcal{K}_n$, 
this proves Step 1 in this case.\\

\noindent Now, assume, for all   infinite 
sequences $(y_n)_n$ in $\mathcal{M}$, there exists $n \not= n'$
such that $\gamma_{y_n} \cap \gamma_{y_{n'}} \not= \emptyset$. Clearly,
this case occurs if $\mathcal{M}$ is finite.
Then, choose a finite set $\{y_1, \cdots, y_k \} \subset \mathcal{M}$ \
such that  
$\gamma_{y_n} \cap \gamma_{y_{n'}} = \emptyset$ if $n \not= n'$ and assume
that $k$ is maximal. This means that for all $y \in \mathcal{M}$, there exists 
$n \in \{1,...,k \}$ such that $\gamma_y  
\cap \gamma_{y_n} \not= \emptyset$ (otherwise the set $\{y,y_1, \cdots, y_k \}$ would have
the above property and $k$ would not be maximal). 
Since $|\gamma_y|=L$, we have:
$\gamma_y \subset \gamma_{y_n}'$ where $\gamma_{y_n}'$ is the trajectory of 
length $3L+2$ and slope $\theta$ 
such that $\gamma_{y_n}'(L+1)= \gamma_{y_n}(0)$. Note that instead of  
$\gamma_{y_n}'$, one could have considered trajectories of length $3L$. We
need to consider trajectories strictly greater than $3L$  to get
relation  
(\ref{3L+2}) below.  For simplicity, we choose $3L+2$.   
Using (\ref{r6}), we see that

\begin{eqnarray} \label{r8}
\tilde{\mathcal{M}} \subset \bigcup_{n=1}^{k} \gamma_{y_n}'.
\end{eqnarray}
Let $n$ be fixed. Let 
$$A_n=\left\{ t \in [0,t_n] \left| \gamma^n(t) \in \bigcup_{m=1}^{k} 
\gamma_{y_m}' \right. \right\}$$ 
The set $ \bigcup_{m=1}^{k} 
\gamma_{y_m}'$
is a finite union of 
disjoint trajectories of 
finite length and slope $\theta$. Moreover, $\theta_n 
\not= \theta$ (see (\ref{r1'})). Hence, $A_n$ is a discrete set of points. 
Since $ \bigcup_{m=1}^{k} 
\gamma_{y_m}'$ is closed, $A_n$ is closed too. This shows 
that $A_n$ is finite. Let now $t_1, t_2 \in A_n$ with $t_1 < t_2$.
As one can check, there exists $\alpha >0$ independent of $n$
such that if $t \leq
\alpha {|\theta_n-\theta|}^{-1}$, $\gamma^n(t_1+t) \not\in 
 \bigcup_{m=1}^{k} 
\gamma_{y_m}'$. Indeed, $\bigcup_{m=1}^{k} 
\gamma_{y_m}'$ is a finite union of segments of slope $\theta$. Then, it is 
clear that the 
trajectory $\gamma_n$, whose slope is $\theta_n$, meets two different 
segments of  $\bigcup_{m=1}^{k} 
\gamma_{y_m}'$ in a 
time proportional to the difference $|\theta_n-\theta|$ of the slopes.   Hence 
$t_2-t_1 \geq \alpha {|\theta_n-\theta|}^{-1}$. This proves that

\begin{eqnarray} \label{r9}
\# A_n \leq \frac{t_n}{\alpha}  |\theta_n-\theta|.
\end{eqnarray} 
Moreover, 
$$m(\gamma^n_{t_n})= \frac{1}{t_n} \int_0^{t_n} \mathcal{X}_{\omega}( 
\gamma^n_{t_n}(t)) dt
=\frac{L}{t_n} \sum_{l=1}^{\frac{t_n}{L}-1} \frac{1}{L}
\int_{lL}^{(l+1)L}  \mathcal{X}_{\omega}( 
\gamma^n_{t_n}(t)) dt=\frac{L}{t_n} \sum_{l=1}^{\frac{t_n}{L}-1} 
m(\gamma^{l,n}).$$
From the definition of $K_n$, we write that for 
$l \in \mathcal{K}_n$ we have $m(\gamma^{l,n}) \geq 0$ and we know that for 
$l \in \{ 0,...,\frac{t_n}{L}-1\} \backslash \mathcal{K}_n$, we have 
$m(\gamma^{l,n}) \geq  g(\om)+ \epsilon$. We obtain that 
$$m(\gamma^n_{t_n}) \geq \frac{L}{t_n}
\#\left(\{ 0,...,\frac{t_n}{L}-1\} \backslash \mathcal{K}_n
\right) (g(\om)+ \epsilon)= \frac{L}{t_n}\left(\frac{t_n}{L}-
\# \mathcal{K}_n \right) (g(\om)+ \epsilon). $$
Since $\lim_n m(\gamma^n_{t_n})=g(\om)$, t 
$$\# \mathcal{K}_n  \geq c t_n$$
where $c>0$ is independent of $n$. It follows from  
(\ref{r9}) that we can choose $k_n \in \mathcal{K}_n$
such that $[k_nL,(k_n+1)L] \cap A_n = \emptyset$.  
Hence, $\gamma^{k_n,n} \cap     \bigcup_{m=1}^{k} 
\gamma_{y_m}' =\emptyset$. 
We can assume that $\lim_n X_{k_n,n}=y \in \mathcal{M}$. Then,
for all $t \in [0,L]$, we have:
$\lim_n \gamma^{k_n,n}(t) = [y,\theta](t)$. Since $\gamma_{y_m}'$ has
been chosen such that their length is $3L+2$, there exists an $m$ such that 
\begin{eqnarray} \label{3L+2}
[y,\theta]_L \subset \gamma_{y_m}'
\end{eqnarray} with $[y,\theta](0) \not=  
\gamma_{y_m}'(0)$ and  $[y,\theta](L) \not=  
\gamma_{y_m}'(L)$. As one can check, this implies that 
$A_{k_n,n}  \cap  \bigcup_{m=1}^{k} 
\gamma_{y_m}' = \emptyset$ (see figure below).


\begin{figure}[h] 
\centerline{\input{carr4.pstex_t}}
\end{figure}


\noindent By (\ref{r7}), we obtain that

\begin{eqnarray} \label{r10}
A_{k_n,n}  \cap  \tilde{\mathcal{M}}
= \emptyset.
\end{eqnarray}

\noindent Now, set
$$d_n=dist(A_{k_n,n},\gamma_{y_m}')=\inf_{X \in A_{k_n,n}} 
dist(X,\gamma_{y_m}')$$
and 
$$d_n'=\sup_{X \in A_{k_n,n}} dist(X,\gamma_{y_m}'). $$
We have $0< d_n < d_n'$.  Since the sequence of sets $(A_{k_n,n})_n$
converges to $\gamma_y$ (see figure above) and since $\gamma_y \subset
\gamma_{y_m}'$, it is clear that  $\lim_n d'n= \lim_n d_n=0$. 
After passing to a subsequence, 
we can assume that
for all $n$, $d_{n+1}' < d_n$. In other words, the sets  
$A_{k_n,n}$ can be assumed to be all disjoint. 
Since $|\partial \omega|$ is finite, there exists
a infinite number of $A_{k_n,n}$ such that 
$|A_{k_n,n} \cap \partial \omega| \leq \epsilon$. We can 
assume that, for all $n$, $|A_{k_n,n} \cap \partial \omega| \leq \epsilon$.
Since $k_n \in \mathcal{K}_n$, this proves Step 1 in this case. Note that
the sequence $(k_n)_n$ does not necessarily tend to $+\infty$.\\

\noindent Let now $\bar{\gamma}^n= \gamma^{k_n,n}$. Keeping the same
notation as in the definition of $A_{k,n}$. Furthermore, let  
$\bar{\gamma}^n_0= \gamma^{k_n,n,0}$. This means that
 $\bar{\gamma}^n_0$
is the unique trajectory of slope $\theta$ and length $L$ such that 
$\bar{\gamma}^n_0(0)= \bar{\gamma}^n(0)$. We now prove that

\begin{step}
We have
$$|m( \bar{\gamma}^n)- m( \bar{\gamma}^n_0)| \leq \frac{2}{L}
\epsilon$$
for any $n$ large enough.
\end{step}
The rays $ \bar{\gamma}^n$ and $ \bar{\gamma}^n_0$ can be seen as segments of 
length $L$ in the plane. Let $p_n$ be the projection on $ \bar{\gamma}^n_0$ 
in the direction of the line $\big(
\bar{\gamma}^n(L) \hbox{ } \bar{\gamma}^n_0(L) \big)$. Thales' theorem
implies 
that  for all $t \in [0,L]$, 
$p_n( \bar{\gamma}^n(t))=  \bar{\gamma}^n_0(t)$.

\begin{figure}[h]
\centerline{\input{carr2.pstex_t}}
\end{figure}



\noindent Now let $t \in [0,L]$. Assume that 
$\mathcal{X}_{\omega}(\bar{\gamma}^n_0(t)) \not= 
\mathcal{X}_{\omega}(\bar{\gamma}^n_0(t))$. Then, there exists a point $X$ of
$\partial \omega$ between $\bar{\gamma}^n(t)$ and $\bar{\gamma}^n_0(t)$. 
Thus,  
$X \in [\bar{\gamma}^n(t),\bar{\gamma}^n_0(t)] \bigcap \partial \omega$.
This shows that $p_n(X) \in p_n(\partial \omega \bigcap 
A_{k_n,n})$. Indeed, remark that, for all $t \in [0,L]$, 
$$  [\bar{\gamma}^n(t),\bar{\gamma}^n_0(t)] \subset A_{k_n,n}.$$
Since $\bar{\gamma}^n_0$ is parametrized by the length, 
the segments $[0,L]$ and $\bar{\gamma}^n_0$ can be identified. We obtain that

\begin{eqnarray} \label{r11}
\big\{t \in [0,L] \hbox{ s.t. } \mathcal{X}_{\omega}(\bar{\gamma}^n(t)) \not= 
\mathcal{X}_{\omega}(\bar{\gamma}^n_0(t)) \big\} \subset 
p_n( \partial \omega \bigcap 
A_{k_n,n}). 
\end{eqnarray}

\noindent
Now, write that 
$$|m( \bar{\gamma}^n)- m( \bar{\gamma}^n_0)| \leq  \frac{1}{L} 
\int_0^L |\mathcal{X}_{\omega}(\bar{\gamma}^n(t))-
\mathcal{X}_{\omega}(\bar{\gamma}^n_0(t))| dt.$$
By (\ref{r11}), this gives 
$$|m( \bar{\gamma}^n)- m( \bar{\gamma}^n_0)| \leq 
 \frac{1}{L}  |p_n( \partial \omega \bigcap 
A_{k_n,n})|.$$
Remember  that $p_n$ is a projection whose angle tends to $\frac{\pi}{2}$. 
Hence, for any piece of curve $C$, 
$|p_n(C)| \leq 2 |C|$ (the constant $2$ could be replaced
here by a constant $c_n$ which goes to $1$ with $n$). It follows that 
$$|m( \bar{\gamma}^n)- m( \bar{\gamma}^n_0)| \leq  
\frac{2}{L}| \partial \omega \bigcap 
A_{k_n,n}|.$$
Step 2 is then a direct consequence of Step 1.

\begin{step}
Conclusion
\end{step}
The trajectory $\bar{\gamma}^n_0$ has length $L$ and slope $\theta$. The 
definition of $L$ (see (\ref{r5}))  then implies that 
$$m( \bar{\gamma}^n_0) \geq g'(\om)- \epsilon.$$
By Step 1, 
$$m( \bar{\gamma}^n) \leq g(\om) + \epsilon.$$
Step 2 then shows that 
$$g'(\om)- \epsilon \leq   g(\om) + \epsilon + 
\frac{2}{L} \epsilon.$$
Since $\epsilon$ is arbitrary, we obtain (\ref{r2}). Together with 
(\ref{r1}), the proof of the theorem is complete.




\end{document}